\documentclass[twoside,reqno,10pt]{amsart}

\newtheorem{theorem}{Theorem}
\newtheorem{lemma}{Lemma}
\newtheorem{corollary}{Corollary}
\newtheorem{proposition}{Proposition}

\newtheorem{definition}{Definition}

\newtheorem{example}{Example}

	\usepackage{amsfonts, amsmath, amssymb,latexsym}

	\usepackage{url}

	\usepackage{tikz-cd}
	\usetikzlibrary{matrix, calc, arrows}
	
	\usepackage{color}
	\definecolor{labelcolor}{RGB}{100,0,0}
	\definecolor{outputcolor}{RGB}{0,0,100}
	
	\usepackage{listings}	
	\definecolor{dkgreen}{rgb}{0,0.6,0}
	\definecolor{gray}{rgb}{0.5,0.5,0.5}
	\definecolor{mauve}{rgb}{0.20 , 0.40, 1.0}
	
	\lstdefinelanguage{Maxima}
	{morekeywords={allbut, block, break, buildq, do, else, elseif, error, errcatch, for, go, if, in, is, local, new, step, then, thru, unless, while, return, true, false, tellsimp, tellsimpafter},
		sensitive=true,
		comment=[s]{/*}{*/},
		morestring=[b]"
	}

\newcommand{\symb}[1]{{\tt #1}}

\newcommand {\grpart}[1]{\ensuremath{\left\langle #1 \right\rangle  }}

\newcommand {\GA}[1] {\ensuremath{\mathbb{G}^{#1} }}
\newcommand {\Cl}[1] {\ensuremath{ C\ell_{#1} } }

\newcommand {\mcommand}[2]{
	{\color{labelcolor}\begin{verbatim}
		(\%i #1) #2
		\end{verbatim}
	}\vskip -1em
}

\begin{document}
 \title[Symbolic Algorithm for Clifford Algebras]{A Symbolic Algorithm for Computation of  Non-degenerate Clifford Algebra Matrix Representations}

\author{Dimiter Prodanov}

\address{
	Neuroelectronics Research Flanders,
	IMEC, Leuven, Belgium \\
	ITSDP, IICT,
	Bulgarian Academy of Science,
	Sofia, Bulgaria
}

\email{dimiterpp@gmail.com; \, dimiter.prodanov@imec.be }

\thanks{The present work is funded  in part by the European Union's Horizon Europe program under grant agreement VIBraTE, grant agreement 101086815.}%

 


\begin{abstract}
Modern advances in general-purpose computer algebra systems offer solutions to a variety of problems, 
which in the past required substantial time investments by trained mathematicians.
An excellent example of such development are the Clifford algebras.
The main objective of the paper is to demonstrate an utterly algorithmic construction of a Clifford algebra matrix algebra representation of a  non-degenerate signature (p, q).
While this is not the most economical way of implementation, it offers a transparent mechanism of translation between a Clifford algebra and its faithful real-valued matrix representation and can be used for automated proof checking.  
This representation is used to derive an algorithm for the computation of an arbitrary multivector inverse as a proof certificate. 
The proposed algorithm is a mapping of the Faddeev--LeVerrier--Souriau algorithm for computation of the characteristic polynomial of matrices.
\end{abstract}

\maketitle


\section{Introduction}
\label{sec:intro}

Modern advances of general-purpose computer algebra systems offer solutions to a variety of problems, 
which in the past required substantial time investments by trained mathematicians.
A good example for such development are the Clifford algebras. 
Recent years have seen renewed interest in Clifford algebra platforms.
Some of these packages have been featured in peer-reviewed publications, 
notably, for Maple \cite{Ablamowicz2002}, Matlab \cite{Sangwine2016}, Mathematica \cite{Aragon-Camarasa2015} and Maxima \cite{Prodanov2016a}.
There is also an inflow of new packages, such as 
\textit{Ganja.js} for JavaScript \cite{Keninck2020}, \textit{GaLua} for Lua 
\textit{Galgebra} for Python 
\textit{Grassmann} for Julia. 

The development of Clifford algebras is based on the insights of Grasssman, Hamilton and Clifford from the 19\textsuperscript{th} century, 
yet at present these algebras are still an active area of mathematical research. 
Clifford algebras provide the natural generalizations of complex, dual  and split complex (or hyperbolic) numbers into the concept of \textit{Clifford numbers}.
Since Clifford algebras are associative, they can be represented by  matrices. 
Examples can be given by the Dirac's $\gamma$ matrices, representing Space Time Algebra, and Pauli's $\sigma$ matrices representing the three-dimensional (3D) algebra of the Euclidean space, that is the Geometric algebra of Clifford. In many applications, the Geometric algebra is referred to as the Algebra of Physical Space.  
The matrix representations have been studied before, resulting in somewhat complicated algorithmic constructions \cite{Brihaye1992,Hile1990,Okubo1991,Okubo1991a,  Okubo1995}.
Recently, the question has been also addressed in \cite{Calvet2017, Lee2013, Song2014}.

The demonstrated application of the presented matrix representation algorithm is the computation of an inverse of a general multivector expression. 
This was achieved in two ways -- 
either by going through an intermediate matrix representation and translating back the obtained result into a Clifford symbolical form or 
by mapping the Faddeev--LeVerrier--Souriau algorithm directly in the algebra itself. 
The paper was presented in preliminary form at the ICCA12 Conference, hosted virtually between 3--7 Aug 2020, at Hefei, China.

Results presented in the paper are derived using an elementary construction of Clifford algebra \cite{Macdonald2002}, 
which is suitable for direct implementation in a computer algebra system supporting symbolical transformations of expressions. 
From design perspective it was the preferred choice in the \symb{clifford} package 
for the computer algebra system Maxima \cite{Prodanov2016,  Prodanov2017, Prodanov2016a}.

The paper is organized as follows.
Sec. \ref{sec:prelim} introduces the notations and conventions.
Sec. \ref{sec:constr} presents a construction of the Clifford algebra, based on \cite{Macdonald2002}.
Sec. \ref{sec:index} introduces the sparse representation of the Clifford algebra.
Sec. \ref{sec:mult} studies the structure of the multiplication table of the algebra.
Sec. \ref{sec:main2} demonstrates the canonical matrix algebra representation theorem.
Sec. \ref{sec:examples} provides examples of low-dimensional matrix algebra  representations.
Sec. \ref{sec:appli} demonstrates the  Faddeev--LeVerrier--Souriau algorithm.
The Appendices demonstrate supporting results. 

\section{Notation and Preliminaries}\label{sec:prelim}

$\Cl{n}$ will mean a Clifford algebra of order \textit{n} but with unspecified signature. 
Algebra generators will be indexed by Latin letters. 
Multi-indices will be denoted with capital letters.
The operation of taking k-grade part of an expression will be denoted by $\left\langle . \right\rangle_k$. 
The notation $(a ; b) \ | \ \mathrm{pred}$ is interpreted as "if the predicate \textit{pred} is true then \textit{a} else \textit{b}".

Matrices will be indicated with bold capital letters, while matrix entries will be indicated by lowercase letters.

\begin{definition}
	The \textbf{extended basis} set of the algebra will be defined as the ordered power set  
	\[
	\mathbf{B} := \left\lbrace P (E), {\prec} \right\rbrace 
	\]	 
	of all generators and their irreducible products. 
\end{definition}
\begin{definition}
	\label{def:mvector}
	A multivector of the Clifford algebra is a linear combination of elements over the $2^n$-dimensional vector space  spanned by $\mathbf{B}$.
\end{definition}

\begin{definition}[Scalar product]\label{def:scprod}
	The \textit{scalar product} of the blades $A$ and $B$ is defined as
	\[
	A \ast B :=  \left\langle A \, B \right\rangle_{0}
	\] 
	and extended by linearity to the entire algebra.
\end{definition}

\begin{definition}[Scalar product table]\label{def:sctab}
	Define the scalar product table
	\[
	\mathbf{G}:= \left\lbrace g_{ \mu \nu} \, e_M \ast e_N \ | M \prec N  \right\rbrace, \quad g_{\mu } \in \{-1,0, 1\}
	\]
\end{definition}

\begin{definition}[Set difference]\label{sec:setdiff} 
	For the sets $A$ and $B$ define the symmetric difference as the operation
	\[
	A  \triangle
	B : = \left\lbrace x: x  \notin A \cap B, \  x  \in A \cup B \right\rbrace 
	\] 
\end{definition}

\section{Elementary construction of the algebra \Cl{p,q,r}}
\label{sec:constr}
 Clifford algebras can be defined in several mathematically equivalent ways.
 However, from computer science perspective the preferred definition should be easy to implement with the tools of the   computer algebra system of choice. 
 This section gives a concise self-contained construction of the Clifford algebra based on Macdonald \cite{Macdonald2002, Prodanov2017}.
 The construction is repeated here for consistency of the presentation. 
 The algebra is constructed by a set of generators, which can be identified with orthonormal vectors.

In order to construct the algebra we need some preliminary conventions. 
We assume that there are
\begin{itemize}
	\item a fixed atomic generator symbol $\mathbf{e}$ 
	\item a set of $n$ abstract \textbf{generators}    $ E:= \left\lbrace e_1 \ldots e_n \right\rbrace $ of the algebra produced by the action of the indicial map
	$ \iota_e: k \mapsto e_k $ for $1 \leq k \leq n$.  
	\item a total order  $ \prec$ over $E$. 
\end{itemize}

\begin{definition}[Clifford algebra]
	\label{def:clialgebra}
 	The construction of the algebra proceeds in the following steps:
	\begin{enumerate}

		\item Define a vector space $V \left(E, \mathbb{R} \right) $ over   
		the set of generators serving as basis with axioms
		\begin{itemize}
			\item Commutativity of vector addition: $a + b = b + a$. 
		    \item Associativity of vector addition:
			$ (a +b ) +c = a + (b+c) $
			\item Existence of additive unity $0$:
		    $ 0 + a = a + 0 =a $.
			\item For every vector, there exists an additive inverse:
		    $a+(-a)=0$.
		    \item Associativity of scalar multiplication:
			$ \lambda (\mu a) = (\lambda \mu) a $
			\item Compatibility with scalar multiplication:
			$a \lambda =  \lambda a $.
		    \item Distributivity of scalar addition:
		    $ (\lambda+ \mu) a= \lambda a+ \mu a$. 	
		    \item Distributivity of vector addition:
		    $ \lambda (a + b)= \lambda a+ \lambda b$.
		    \item Scalar multiplication identity:
		    $ 1_{\mathbb{R}}\; a = a $.		    
		\end{itemize}    
	for vectors $a,b,c$ and scalars $\lambda, \mu$. 
	\item Adjoin an associative algebra $ \mathbb{G} (E, \mathbb{R} )$ over  $V \left(E, \mathbb{R} \right) $ using the \textsl{Clifford product} or \textsl{geometric multiplication} operation with axioms
	 \begin{itemize}
		\item Existence of an algebra unity: $1_{\Cl{}} a = a 1_{\Cl{}} =  a $ for a non-scalar $a$.
	 	\item Left distributivity: $a  (b + c) = a  b + a  c$ for arbitrary elements $a, b, c$.
	 	\item Right distributivity: $ (a + b)  c = a b + a  c$ for arbitrary elements $a, b, c$.
	 	\item Associativity: $a (b c)=  (a b) c$ for arbitrary elements $a, b, c$.
	 	\item Compatibility with scalars: $(\lambda a)  (\mu b) = (\lambda \mu ) (a  b)$ 
	 for scalars $\lambda, \mu$ and non-scalar elements $a, b$.
	 \end{itemize}
	 
	 \item Finally, assert
		\begin{description}
			\item[Closure Axiom]
			For $k$  generators $e_1, \ldots, e_k \in E$ and  scalars $\lambda$, $\mu$ 
			the multivector belongs to the algebra:
			\begin{equation*}
			\lambda	+ \mu e_1 \ldots e_k \in \GA{n}   \tag{C}
			\end{equation*}

			\item[Reduction Axiom]  
			For all generators
			\begin{equation}\label{eq:square}
			e_k e_k= \sigma_{k} 1_{\Cl{}}, \ \sigma_{k} \in \left\lbrace 1, -1, 0 \right\rbrace \tag{R}
			\end{equation}
			where $\sigma_k   $ are real  or complex scalars.  
			The signature set is $S:= \{\sigma_k \}$.

			\item[Anti-Commutativity Axiom]
			For every two basis vectors, such that $ e_i \prec e_j$
			\begin{equation}\label{eq:anticomm}
			e_i e_j=  - e_j e_i \tag{A-C}
			\end{equation}

		\end{description}
	\end{enumerate}
We denote the algebra satisfying axioms 1-3 as Clifford algebra  \Cl{p,q,r}, where $p$ elements of the orthonormal basis square to $1$, $q$ elements square to $-1$ and $r$ (degenerate) elements square to $0$.
\end{definition}

It should be remarked that the closure and compatibility Axioms are not included in the original construction of Macdonald.
The algebra construction can be extended to other fields with characteristic zero. 
So-presented construction can be carried out without modifications in Computer Algebra systems like \textit{Maxima}
by defining proper simplification rules for the geometric product operation \cite{Prodanov2016,Prodanov2016a}.

The definition does not specify a particular ordering but merely asserts that an ordering must be used.


\subsection{Consistency of the extension}
\label{sec:cons}

Here we present results demonstrating the consistency of the Clifford algebra definition. 
 \begin{proposition}
	The scalar and algebra units coincide:
	\[
	1_\mathbb{R}= 1_{\Cl{}} \equiv 1
	\] 
\end{proposition}
\begin{proof}
	Suppose that $e_0  \in V (E) $ is the unit of the algebra. Then by the anti-commutativity: 
	$e_0 e_k + e_k e_0 =0  \leftrightarrow e_k +e_k=0$ which is a contradiction.
	Therefore, $e_0 \in \mathbb{K} \Rightarrow e_0 =1$ by the uniqueness property.
\end{proof}

\begin{theorem}
	\label{th:vextend}
	$V \left( {E}  \right) $ extends over the (ordered) power set of $E$.
	\[
	P (E)  := \{ e_1, \dots, e_n, e_1 e_2, e_1 e_3, \ldots ,e_1 e_2...e_n\} 
	\]
	That is $V \left( P(E)  \right)$  is a vector space and we have the inclusion
	\[
	V \left( {E}  \right)  \subset V \left( P(E)  \right)
	\]
\end{theorem}
\begin{proof}
	The proof follows from the Closure Axiom.
	
	Multiplicative properties:
	\begin{enumerate}
		\item Scalar commutativity $ \mu e_1 \ldots e_k = e_1 \mu \ldots e_k = \ldots =  e_1 \ldots e_k \mu $
	
		\item In particular $ 1 e_1 \ldots e_k =   e_1 \ldots e_k 1 = e_1 \ldots e_k $.	
	
		\item $0 , e_1 \ldots e_k = 0 \, e_2 \ldots e_k =\ldots = 0 \, e_k =0$.
	
		\item Scalar compatibility: $ \lambda ( \mu  e_1 \ldots e_k) = \lambda ( \mu  e_1) (e_2 \ldots e_k) = (\lambda \mu ) e_1 \ldots e_k  $
	\end{enumerate}
	And we extend by linearity and the closure over the whole algebra. 
	
	Additive properties:
	\begin{enumerate}
	 \item Universality of zero: Let $L = 0 + e_1 \ldots e_k$. Right-multiply by $e_k$.
	 If $\sigma_k=0$ then the results follows trivially. 
	 
	 So let $\sigma_k \neq 0$.
	 Then $L e_k = 0 e_k + e_1 \ldots e_{k-1} \sigma_k=  \sigma_k e_1 \ldots e_{k-1}  $.
	 Right-multiply by $e_k$. Then
	   $ \sigma_k L =   \sigma_k  e_1 \ldots e_k   \Leftrightarrow L= e_1 \ldots e_k $.
	   
	 \item Scalar distributivity: $ (c+d) e_1 \ldots e_k = ((c+d) e_1 ) e_2 \ldots e_k = (c e_1 +d e_1) e_2 \ldots e_k =
	 c e_1  \ldots e_k  +d e_1  \ldots e_k$.
	 
	 \item  In particular, there is an additive inverse element 
	 $ e_1  \ldots e_k -  e_1  \ldots e_k =0 $.
	 \item Commutativity of addition:
	 Let $M= e_a + e_b - e_b - e_a  $. 
	 Then for all possible orders of evaluation $M=0$.
	 Therefore, by linearity $a + b = b + a$. 
	
	\item Associativity of addition:
	Let $M= (e_a + e_b) +  e_c - e_a - (e_b + e_c)$. By commutativity of addition and closure
	$M=(e_a + e_b) 1 - (e_b + e_c) 1 -  e_c - e_a  $.
	Then for all possible orders of evaluation of the first two summands
	$M =e_a - e_c + e_c - e_a =0$.
	Therefore, by linearity $(a + b) +c = a+ (b + c)$. 
    \item Distributivity of addition follows directly from the distributivity axioms and the closure axioms
    \[
    \lambda (e_a +e_b) = \lambda e_a  + \lambda e_b
    \]
	\end{enumerate}
	And we extend by linearity and the closure over the entire algebra. 
\end{proof}

Based on this extension principle we can promote the indical map to a full index set isomorphism
\begin{equation}\label{eq:indexisom}
j : e_a  e_b  \ldots e_m  \mapsto e_{J}, \quad J =\left\lbrace a, b, \ldots, m \right\rbrace 
\end{equation}
and use  multi-index notation implicitly wherever appropriate. 
For consistency, we also extend the ordering to the power set
$
  \prec_{ P(E)} 
$, however no distinct symbol will be used for this.

\begin{theorem}[Well-posedness]\label{th:wellpose}
	The Clifford algebra construction specified in Def. \ref{def:clialgebra} is well-defined. 
\end{theorem}
\begin{proof}
	$V (E)$ is well-defined. Moreover, the vector space extension is also well-defined.
	The compatibility Axioms imply the inclusion
	$  V (E, \mathbb{R} ) \subset  \mathbb{G} \left(E, \mathbb{R} \right) $.
 
	Therefore, we only need to check the additional Algebra axioms.  	
	The Reduction Axiom agrees with the Closure Axiom trivially.
	The Anti-commutativity Axiom agrees with the Closure Axiom trivially.
	Further, by restriction to ordered pairs the  
	Anti-commutativity Axiom does not apply when the Reduction Axiom applies. 
	Therefore, these two axioms are logically independent. 
	Therefore, restricted algebra construction over $V(E)$ is well-defined. 
	The algebra has a maximal element by Prop. \ref{prop:pscalar}.
	The algebra extension is well-defined since the existence of \Cl{p+q+r=n} implies the existence of \Cl{p+q+r=n-1} by the Closure Axiom.
	By reduction, for $n=1$ the three possible cases for the sign of the $e_1^2$  represent the  double, complex or dual numbers, respectively. 
	Therefore, by induction \Cl{p+q+r=n} exists.
\end{proof}

\begin{theorem}[Algebra equivalence]
	\label{th:aleq}
	Two algebras with the same numbers of $p$, $q$ and $r$ generators are order-isomorphic.
\end{theorem}
\begin{proof}
	If the elements are ordered identically then the statement is trivial. Let's assume that the elements are ordered differently.
	But then there is a permutation putting them into canonical order. 
	Therefore, the algebras are isomorphic. Therefore, we can identify an order $\prec^\prime$ with the second permutation. 
\end{proof}

\begin{corollary}[S-Law of inertia]
	Denote by \textit{S} the signature set of the algebra \Cl{}.
	For two isomorphic algebras ${\Cl{}}^\prime $ and ${\Cl{}}^{\prime\prime} $ ( ${\Cl{}}^\prime \cong {\Cl{}}^{\prime\prime}$ ) there is an invertible map
	\[
	P : S_1 \mapsto S_2 
	\]
	Conversely, if there is permutation 
	\[
	P : S_1 \mapsto S_2 
	\]
	then   ${\Cl{}}^{\prime}  \cong {\Cl{}}^{\prime\prime}$.  
\end{corollary}
\begin{proof}
	The forward statement follows from Th. \ref{th:aleq}.
	Converse case:
	If the dimensions of $S_1$ and $S_2$ are equal and the numbers of $p$, $q$ and $r$ generators are equal
	then there is a permutation
	\[
	P : S_1 \mapsto S_2 
	\] 
	However, it can be noted that this is the same permutation which gives 
	\[
	P : \mathbf{B}_1 \mapsto \mathbf{B}_2 
	\]
	and hence by Th. \ref{th:aleq} $ {\Cl{}}^{\prime} \cong {\Cl{}}^{\prime\prime}$.
\end{proof}
\begin{definition}[Canonical real algebra]
	Define the canonical ordering as the nested lexicographical order $\varrho$, such that  
	$i< j \Longrightarrow e_i \prec e_j$ and extend it over $P(E)$ as:
	\[
	e_1 \prec e_2 \prec \underbrace{e_{1} e_2}_{e_{1 2}} \prec \ldots \prec \underbrace{ e_1 \ldots e_n}_{e_{N}}
	\]
	In addition assume that the first $p$ elements square to 1, the next $q$ elements square to -1 and the last $r$ elements square to 0.
	Then the algebra 
	$
	\Cl{p,q,r} $ specified by the structure $\left\lbrace e, \varrho,  \mathbb{R} \right\rbrace 
	$ is the canonical Clifford algebra.
\end{definition}
\begin{corollary}
	The canonical algebra $\Cl{p,q,r} $ defines an equivalence class with regard to permutations.
\end{corollary}

From computational perspective, last two corollaries are important because they indicate that an algorithm cane be implemented only for one Clifford algebra of a certain grade. 
They also imply that the order in which operations are composed can substantially simplify computations.

\section{Indicial or sparse representation}
\label{sec:index}

\begin{definition}[Indicial map]		
	Define the indicial map $\iota_e$, acting on symbols by concatenation (i.e. of a set), such that
	\begin{flalign*}
	\iota_e : & \oplus_n \mathbb{N} \mapsto P (E)  := \{ e_1, \dots, e_n, e_1 e_2, e_1 e_3, \ldots ,e_1 e_2...e_n\}  \\
	\iota_e : & g   \mapsto e_g
	\end{flalign*}
where g is set-valued 
	and assert the convention $ \iota_e : \ \emptyset   \mapsto 1_{Cl} $.
\end{definition}


\begin{definition}
	Define the argument map $\mathrm{arg}$ acting on symbol compositions as
	\[
	\mathrm{arg}: \ f(g) \equiv f \circ g \mapsto g
	\]
	Assert $ \mathrm{arg} f = \emptyset $ for the atomic symbol $f$. 
\end{definition}

These definitions allow for stating a very general result about Clifford algebra representations.
\begin{theorem}[Indicial representation]\label{th:indic}
 	For generators $e_s, e_t \in  \Cl{p,q,r}$, such that $ s \neq t$, the following diagram commutes
	
\begin{center}
	\begin{tikzpicture}
	\matrix (m) [matrix of math nodes, row sep=3em,column sep=4em, minimum width=2em, minimum height=2em, nodes={asymmetrical rectangle}]
	{
		e_s &    \{ s \}  \\
		e_s e_t \equiv e_{st} &   \{ s, t \}  \\
	};
	\path[->, font=\scriptsize ]
	([yshift=0.5ex] m-1-1.east)  edge node [above] { $\mathrm{arg}$}   ([yshift=0.5ex] m-1-2.west) ;
	\path[ -> ]
	([yshift= -0.5ex] m-1-2.west)  edge   node [below] { $\iota_e$} ([yshift=-0.5ex] m-1-1.east) ;
	\path[->, font=\scriptsize ]
	(m-1-1)       edge node [right] { $ e_t $}  (m-2-1);
	\path[->, font=\scriptsize ]        
	(m-1-2) edge node [right] { $ \triangle \{ t \}$ } (m-2-2)  ;
	\path[->, font=\scriptsize ]
	([yshift=0.5ex] m-2-1.east)  edge node [above] { $\mathrm{arg}$}   ([yshift=0.5ex] m-2-2.west) ;
	\path[ -> ]
	([yshift= -0.5ex] m-2-2.west)  edge   node [below] { $\iota_e$} ([yshift=-0.5ex] m-2-1.east) ;
	\end{tikzpicture}
\end{center}
	
\end{theorem}
\begin{proof}
	 The right-left $\iota$ action follows from the construction of \Cl{p,q,r}.
	 The left-right argument action is trivial.
	 We observe that $ \mathrm{arg} f = \emptyset $. 
	 Trivially,  $ \{s\}  \triangle  \emptyset = \emptyset  \triangle  \{s\}  = \{s\} $.
	 Let's suppose that  $s=t$.
	 We notice that   $ \{s\}  \triangle  \{t\}= \{t\}  \triangle  \{s\}  = \emptyset $.
	 Let's suppose that  $s \neq t$.
	 We notice that   $ \{s\}  \triangle  \{t\}=      \{s,  t\}  $ and
	 $ \{t\}  \triangle  \{s\}=      \{t,  s\}  $.
\end{proof}
This theorem can be used for the reduction of Clifford products.


\subsection{Simplification of Clifford products}\label{sec:simp}

This section presents how to compute the \textbf{canonical simplified form} of blades.
The non-simplified arbitrary Clifford product of multivectors will be referred to simply as a Clifford multinomial.
The main lemma is demonstrated in \cite{Macdonald2002} and is repeated here for convenience:
\begin{lemma}[Permutation equivalence]
	\label{th:permeq}
	Let $ M= e_{k_1} \ldots e_{k_i} $ be a Clifford multinomial, where the indices are not necessarily different. 
	Then 
	\[
	M = s \, P_\rho \left\lbrace  e_{k_1} \ldots e_{k_i} \right\rbrace 
	\]
	where $s =\pm 1$ is the sign of permutation of $M$ and 	$ P_\rho \left\lbrace  e_{k_1} \ldots e_{k_i} \right\rbrace  \mapsto  e_{k_\alpha} \ldots e_{k_\omega}$ is the  product permutation according to the canonical ordering.
\end{lemma}
\begin{proof}
	The proof follows directly from the anti-commutativity of Clifford multiplication for any two generator elements \ref{eq:anticomm}, observing that 
	the sign of a permutation of $S$ can be defined from its decomposition into the product of transpositions as
	$sgn(B) = (-1)^m$, where $m$ is the number of transpositions in the decomposition. 
\end{proof}
Further, we can define a \textit{simplified form} according to the action of the Reduction Axiom. 
\begin{theorem}[Reduced blade form]
	\label{th:simpf}
	Let $  e_{k_1} \ldots e_{k_i} $ be an arbitrary Clifford multinomial, where the $i$ generators in \Cl{p,q,r} are not necessarily different.
	Let $s_P$ be the sign of the permutation 	$ P_\rho \left\lbrace  e_{k_1} \ldots e_{k_i} \right\rbrace $.
	We say that simplification induces a reduced or simplified form $\mathcal{R}$ with regard to permutations and product evaluations
	such that
	\begin{flalign*}
	  e_{k_1} \ldots e_{k_i} = & s_P \  \mathcal{R} \left( e_{k_\alpha} \ldots e_{k_\omega} \right)
	=  s_P \ e_{k_\alpha} \mathcal{R} \left( e_{k_\beta} \ldots e_{k_\omega}   \right) =  \\
	& 	s_P \ \mathcal{R} \left( e_{k_\alpha} e_{k_\beta} \right) \mathcal{R} \left( e_{k_\gamma}  \ldots e_{k_\omega}   \right)
	\end{flalign*}
	Let $2 p \leq i$ generators square as $e_j^2=1$. 
	Then its simplified form is  
	\[
	\mathcal{R} \left( e_{k_1} \ldots e_{k_i} \right) =  s_P \,   e_{M} , \quad M = \left\lbrace k_1 \triangle k_2 \triangle \ldots  k_{i - 2 p} \right\rbrace 
	\]
	with all $k$-indices different. 	
	Let $2 q \leq i$ generators square as $e_j^2=-1$. 
	Then its simplified form is    
	\[
	 \mathcal{R} \left(e_{k_1} \ldots e_{k_i}\right)  =  (-1)^q s_P \  e_{M} , \quad M = \left\lbrace k_1 \triangle k_2 \triangle \ldots  k_{i - 2 q} \right\rbrace 
	\]
	with all $k$-indices different. 
	Let at lest $2 $ generators square as $e_j^2=0$. 
	Then its simplified form is  
	\[
	\mathcal{R} \left(	e_{k_1} \ldots e_{k_i} \right) =0
	\]
\end{theorem}
\begin{proof}
	The product is transformed according to Lemma \ref{th:permeq}.
	Then we use Th \ref{th:indic} to compute the indices. 
	By the Reduction axiom (\ref{eq:square}) the elements of equal indices are removed from the final index list. 
	This process induces a factor of $(-1)^q$ in the final result. 
	Hence, for
	$2 q $ elements of index $a$ $ \left\lbrace k_a \right\rbrace \triangle \left\lbrace k_a \right\rbrace = \emptyset $. 
	Further
	$ \emptyset \triangle \left\lbrace k_i \right\rbrace = \left\lbrace k_i \right\rbrace  \triangle \emptyset = \left\lbrace k_i \right\rbrace $. 
	The final list can be translated back into a blade form by the map $\iota_e$ by Th. \ref{th:indic}.
	The associativity of $\mathcal{R}$ follows from the associativity of the Clifford product.
\end{proof}
The form equivalence is used for simplification of expressions in the Clifford Maxima package.


\section{Multiplication tables of  Clifford algebras}
\label{sec:mult}

\begin{definition}[Full Matrix multiplication table]
	Consider the extended basis $\mathbf{B}$. Define the multiplication table matrix as the mapping of the basis 
	into the set of square matrices valued in \Cl{p+q+r}: 
	\begin{flalign*}
	\mu & :   \mathbf{B} \times \mathbf{B}  \mapsto \mathbf{Mat} ({2^{n}} \times  {2^{n}} ), \ n=p+q+r	\\
	\mu &  (\mathbf{B})  =  \mathbf{M}_{\Cl{p,q,r}}  
	\end{flalign*}	 
	with matrix consisting of the ordered ans simplified (i.e. $\mathcal{R}$ form) product entries using the multi-index notation
	\[
	\mathbf{M} := \left\lbrace m_{ \mu \nu} \, \mathcal{R} (e_M e_N) \ | M \prec N  \right\rbrace, \ \ m_{ \mu \nu}=\left\lbrace -1, 0, 1 \right\rbrace 
	\]
\end{definition}
Since $\mathbf{M}$ is a square matrix we will not make a distinction for  its dimensions, so $dim(\mathbf{M})$ will mean
$ 2^n \times 2^n $ or $2^n $ depending on the context.

\begin{example}
	Consider $\Cl{2,0,0}$. The full matrix multiplication table is
	\[
	\mathbf{M}= \left( \begin{array}{c|cccc}
	1 & {e_1} & {e_2} & {e_1} {e_2}\\
	\hline
	{e_1} & 1 & {e_1} {e_2} & {e_2}\\
	{e_2} & -{e_1} {e_2} & 1 & -{e_1}\\
	{e_1} {e_2} & -{e_2} & {e_1} & -1
	\end{array}\right) 
	\] 
	The lines are drawn for clarity.
\end{example}

\begin{proposition}\label{prop:msparse}
	Consider the multiplication table $\mathbf{M}$.
	All elements $m_{k  j} $ are different for a fixed row $k$.
	All elements $m_{i q }$ are different for a fixed column $q$.
\end{proposition}
\begin{proof}
	The results follow from the simplification Th. \ref{th:simpf}.
	Fix $k$. Then for $e_K, e_J \in \mathbf{B} $
	\[
	e_{K} e_{J} = m_{k  j} e_s, \ \ S = K \triangle J
	\]
	which is a different way of writing the simplification Lemma. 
	Suppose that we have equality for 2 indices $j, j^\prime$. Then
	\[
	K \triangle J^\prime = K \triangle J = S
	\]
	and let $\delta = J \cap J^\prime$. Then 
	\[
	K \triangle  \left( J \cup \delta \right)  = K \triangle J = S \Rightarrow K \triangle \delta = S \Rightarrow \delta = \emptyset
	\]
	Therefore, $j = j^\prime$. 
	By symmetry, the same reasoning applies to a fixed column $q$.
\end{proof}

\begin{lemma}[Multiplication Matrix Structure]\label{th:struct}
	For the  multi-indices $S \prec T$, such that $ S \cap T =   \emptyset   $, 
	the following implications for the   elements of $\mathbf{M}$ hold:
	\begin{center}
		
		\begin{tikzpicture}
		\matrix (m) [matrix of math nodes,row sep=3em,column sep=4em,minimum width=2em]
		{
			m_{\mu \lambda}\,  e_s &  m_{\mu   \lambda^\prime} e_t & \\
			m_{\lambda \mu} e_s &  m_{\lambda \mu} m_{\mu \lambda^{\prime } } e_{s t} &   m_{\lambda \lambda^{\prime \prime} } e_{s t} \\};
		\path[-> ]
		(m-1-1) edge node [left] { $\exists  $} (m-2-1)
		edge node [above] {$\exists \lambda^\prime > \lambda $} (m-1-2)
		(m-1-2) edge (m-2-2)
		(m-2-1) edge node[above] {$\exists$} (m-2-2)
		(m-2-2) edge node[above]{$\exists \lambda^{\prime \prime}= \lambda^\prime $}  (m-2-3);
		\end{tikzpicture} 
	\end{center}
	so that 
	\[
	m_{ \lambda \lambda^\prime  } \, e_{s t} =  m_{\lambda \mu }  \sigma_\mu m_{\mu \lambda^\prime }\, e_{s t} 
	\]
	for some index $\mu$.
\end{lemma}

\begin{proof}
	By the properties of $\mathbf{M}$ there exists $\lambda^\prime > \lambda$, such that
	$e_{M} e_{L^\prime}   = m_{\mu \lambda^\prime }\,  e_t, \ \  L^\prime \backslash M  = T$
	for $L \prec   L^\prime $.
	Choose $M$, s.d. $L \prec M \prec  L^\prime $.
	Then for
	 $L \prec M \prec  L^\prime $ and $S \prec T$
	\begin{flalign*}
	e_{M} e_{L} & = m_{\mu \lambda}\,  e_s , \ \ {\small L \triangle M } = S
	\Leftrightarrow e_{L} e_{M} = m_{\lambda \mu }\,  e_s \\
	e_{M} e_{L^\prime} & = m_{\mu \lambda^\prime }\,  e_t, \ \ {\small L^\prime \triangle M } = T
	\end{flalign*}
	
	Suppose that $e_s e_t =  e_{s t}, st = S \cup T = S \triangle T$. 
	Multiply together the diagonal nodes in the matrix
	\[
	e_{L} \underbrace{e_{M}  e_{M}}_{\sigma_\mu} e_{L^\prime} = m_{\lambda \mu }  m_{\mu \lambda^\prime }\, e_{s t} 
	\]
	Therefore,
	$ s\in L $ and
	$ t \in L^\prime $.
	We observe that there is at least one element (the algebra unity) with the desired property $\sigma_\mu \neq 0$.
	
	Further, we observe that there exists unique $\lambda^{\prime \prime} $ such that $  m_{\lambda \lambda^{\prime \prime} } e_{s t} $.
	Since $\lambda $ is fixed by Th. \ref{th:simpf} this implies that $ L^{\prime \prime} = L^\prime \Rightarrow \lambda^{\prime \prime} = \lambda^{\prime}$.
	Therefore,  
	\[
	e_{L} e_{L^\prime} = m_{ \lambda \lambda^\prime  } e_{s t},  \ \  L^\prime \triangle L  = \{ s, t \} 
	\]
	which implies the identity
	\[
	m_{ \lambda \lambda^\prime  } \, e_{s t} =  m_{\lambda \mu }  \sigma_\mu m_{\mu \lambda^\prime }\, e_{s t} 
	\]
\end{proof}

\section{Clifford algebra real matrix representation theorem}
\label{sec:main2}

In the present article we will focus on non-degenerate Clifford algebras. 
Therefore, we assume that the product set 
\[
S := \left\lbrace \sigma_i = e_i e_i \ | \ e_i \in \mathbf{E} \right\rbrace   
\]
is valued in the set $\left\lbrace -1, 1 \right\rbrace$ if not stated otherwise.

In order to state the main result we need the following definitions:
\begin{definition}[Coefficient map]\label{def:coefmap}
	Define the linear map acting element-wise
	\[
	C_a : \Cl{} \mapsto \mathbb{R}  
	\]
	by the action
	\[
	C_a:  \begin{cases} a x &\mapsto x, x \in \mathbb{R}, a \in \mathbf{B} \\
	b &\mapsto 0, b \in \mathbf{B}
	\end{cases}
	\]
%
	Define the \textbf{coefficient map} indexed by the multi-index $S$ as
	\[
	C_S : \mathbf{M} \mapsto \mathbf{A}_S
	\]
\end{definition}

\begin{definition}[Canonical matrix map]
	\label{def:matcanon}
	For the  multi-index $S$ define the map
	\[
	\pi : e_{S} \mapsto \mathbf{E}_s= \mathbf{G} \mathbf{A}_s
	\]
	where \textit{s} is the ordinal of $e_{S}$ in the multivector basis $\mathbf{B}$ and $\mathbf{A}_S$ is computed as in Def. \ref{def:coefmap}.
	Further, denote the set of all maps as $\pi= \left\lbrace \pi_s \right\rbrace $ and let $\pi_s \equiv \pi(e_s)$.  
\end{definition}

\begin{proposition}\label{prop:coeflin}
	The  $\pi$-map is linear.
\end{proposition}
The proof follows from the linearity of the coefficient map and matrix multiplication.

\begin{theorem}[Semigroup property]
	Let $e_s$ and $e_t$ are basis elements.
	Then the map $\pi$ acts on \Cl{p,q} according to the following diagram
	
	\begin{center}
		
		\begin{tikzpicture}
			\matrix (m) [matrix of math nodes, row sep=3em,column sep=4em, minimum width=2em, minimum height=2em, nodes={asymmetrical rectangle}]
			{
				e_s &   \mathbf{E}_s   \\
				e_s e_t \equiv e_{st} &  \mathbf{E}_{st} \equiv \mathbf{E}_{s} \mathbf{E}_{t} \\
			};
			\path[->, font=\scriptsize ]
			(m-1-1) edge  node [right] { $e_t$  } (m-2-1)
			edge   node [above] {$\pi $} (m-1-2)
			(m-1-2) edge node [right] { $\mathbf{E}_t$  }(m-2-2)   
			(m-2-1) edge node[above] {$\pi$} (m-2-2);
		\end{tikzpicture}
		
	\end{center}

	The map $\pi$ distributes over the Clifford product:
	\[
	\pi (e_s e_t) = \pi (e_s) \pi(e_t) 
	\]
	The set of all matrices $\mathbf{E}_s$ forms a semigroup.
\end{theorem}
\begin{proof}
	Let 
	\begin{flalign*}
		\pi : \ e_s & \mapsto \mathbf{E}_s \\
		\pi :  \ e_t & \mapsto \mathbf{E}_t \\	
		\pi :  \ e_s e_t & \mapsto \mathbf{E}_{st} \\
	\end{flalign*}
	We specialize the result of Lemma \ref{th:struct}  for $S= \{ s\}$ and $T= \{ t\}$ and observe that
	\[
	m_{ \lambda \lambda^\prime  } \, e_{s t} =  m_{\lambda \mu }  \sigma_\mu m_{\mu \lambda^\prime }\, e_{s t} 
	\]
	for  $ \lambda, \lambda^\prime, \mu  \leq n$
	\[
	\sigma_\lambda m_{ \lambda \lambda^\prime  } =  \sigma_\lambda m_{\lambda \mu }  \sigma_\mu m_{\mu \lambda^\prime } 
	\]
	Therefore,
	\[
	\mathbf{E}_{s t} = \mathbf{E}_s \mathbf{E}_t
	\]
	Moreover, we observe that 
	\[
	\pi (e_s e_t) =\mathbf{E}_{s t}= \mathbf{E}_s \mathbf{E}_t = \pi (e_s) \pi(e_t) 
	\]
	For the semi-group property, consider that 
	since $\pi$ is linear it is invertible. 
	Since $\pi$ distributes over Clifford product its inverse $\pi^{-1}$ distributes over
	matrix multiplication:
	\[
	\pi^{-1} (\mathbf{E}_{s} \mathbf{E}_{t} ) \equiv  \pi^{-1} (\mathbf{E}_{s t}) = e_{st} \equiv e_{s} e_{t} = \pi^{-1} (\mathbf{E}_s) \ \pi^{-1} (\mathbf{E}_t)
	\]
	However, $\Cl{p,q}$ is closed by construction, therefore, the set $\left\lbrace \mathbf{E}\right\rbrace _s$ is closed under	matrix multiplication.
\end{proof}
\begin{corollary}\label{corr:pinn}
	For $\Cl{p,q}$ the set of all matrices $\mathbf{E}_s$ forms a group, which is a representation of the Pin group $ \mathrm{Pin}_{p,q}$.
\end{corollary}
\begin{proof}
We observe that for non degenerate algebras every generator has an inverse element by Props. \ref{prop:vsquare};
therefore for non-degenerate algebras \Cl{p, q}
\[
\mathbf{E}_s^{-1}= \sigma_s \mathbf{E}_s
\]
However, $\sigma$ is valued in $\{ -1, +1\}$.
\end{proof}
\begin{proposition}
	Let $\mathbf{L}:= \left\lbrace l_i \right\rbrace, l_i \in P (E)$ and
	$\mathbf{R}_s$ is the first row of $\mathbf{E}_s$.
	Then 
	$\pi^{-1}: \mathbf{E}_s \mapsto  \mathbf{R}_s \mathbf{L}$.
\end{proposition}
\begin{proof}
	We observe that by the Prop. \ref{th:sparsec} the only non-zero element in the first row of $\mathbf{E}_s $ is 
	$\sigma_1 m_{1 s} = 1$.
	Therefore, 
	$ \mathbf{R}_s \mathbf{L} = e_s $.
\end{proof}

The main result of the section is stated in the theorem below.
\begin{theorem}[Canonical Real Matrix Representation]
Define the map
$g: \mathbf{A} \mapsto \mathbf{G} \mathbf{A}$.
Then
\[
\pi_s = C_s \circ g = g \circ C_s
\]
so that the diagram commutes
\begin{center}
\begin{tikzpicture}
	\matrix (m) [matrix of math nodes, row sep=3em,column sep=4em, minimum width=2em, minimum height=2.5em, nodes={asymmetrical rectangle}]
	{
		\mathbf{M} &   \mathbf{A}_s   \\
		\mathbf{G} \mathbf{M} &  \mathbf{E}_{s} = \mathbf{G} \mathbf{A}_s \\
	};
	\path[->,font=\scriptsize ]
	(m-1-1) edge node [right] { $g$ }  (m-2-1)
	        edge node [above] { $C_s$} (m-1-2)
	(m-1-2) edge node [right] { $g$ }  (m-2-2)   
	(m-2-1) edge node [above] { $C_s$} (m-2-2);
	\path[-> ]
	(m-1-1) edge node [auto] { $\pi_s$} (m-2-2);
	\end{tikzpicture}
\end{center}

Further, $\pi$ is an isomorphism inducing a Clifford algebra representation in the real matrix algebra:
\begin{center}
\begin{tikzpicture}
\matrix (m) [matrix of math nodes, row sep=3em,column sep=4em, minimum width=2em, minimum height=2.5em, nodes={asymmetrical rectangle}]
{
	C\ell_{p,q} (\mathbb{R}) &   C\ell_{p,q} \left[ \mathbf{Mat}_{\mathbb{R}} ({2^{n}} \times  {2^{n}}   ) \right]   \\
};
\path[ -> ]
([yshift=0.5ex] m-1-1.east)  edge node [above] { $\pi$}   ([yshift=0.5ex] m-1-2.west) ;
\path[ -> ]
([yshift= -0.5ex] m-1-2.west)  edge   node [below] { $\pi^{-1}$} ([yshift=-0.5ex] m-1-1.east) ;
\end{tikzpicture}
\end{center}

The $\pi$ map distributes over the Clifford product (homomorphism):
\[
\pi_{st} = \pi_s \pi_t
\]

\end{theorem}
\begin{proof}
	The  $\pi$-map is a linear isomorphism.
	The set $ \left\lbrace \mathbf{E}_s \right\rbrace  $ forms a group, which is a subset of the matrix algebra 
	$ \mathbf{Mat}_{\mathbb{R}} (N \times N), N=2^n $. 
	Let 
    $ \pi (e_s)= \mathbf{E}_s$ and $
	\pi (e_t)= \mathbf{E}_t 
	$.
	It is claimed that
	\begin{enumerate}
		\item $ \mathbf{E}_s \mathbf{E}_t \neq \mathbf{0}$ by the Sparsity Lemma \ref{th:sparsity}.
		\item $ \mathbf{E}_s \mathbf{E}_t =- \mathbf{E}_t \mathbf{E}_s$ by Prop. \ref{prop:antisym}.
		\item $ \mathbf{E}_s \mathbf{E}_s = \sigma_s \mathbf{I}$ by Prop. \ref{prop:vsquare}.
	\end{enumerate}
Therefore, the set $ \left\lbrace \mathbf{E}_s \right\rbrace_{s=1}^{N}  $ is an image of \Cl{p,q}.
\end{proof}

So-constructed canonical matrix representation has an interesting structure.
This can be summarized in the following results.
\begin{proposition}\label{prop:cases}
	Consider the  non-degenerate algebra \Cl{p,q}. Then for any element $e_s$ :
	\begin{itemize}
		\item If $\sigma_s = 1 $ then $ \mathbf{E}_s$ is symmetric.
		\item If $\sigma_s = -1 $ then $ \mathbf{E}_s$ is anti-symmetric.
	\end{itemize}
\end{proposition}
\begin{proof}
	\[
	\mathbf{E^{-1}}_s = \sigma_s \mathbf{E}_s = \mathbf{E}_s^T \Rightarrow  \sigma_s \mathbf{E}_s - \mathbf{E}_s^T = \mathbf{0}
	\]
	Then for $\sigma_s = 1 \Rightarrow   \mathbf{E}_s - \mathbf{E}_s^T  = \mathbf{0}$, therefore $ \mathbf{E}_s$ is symmetric.
	While for
	$\sigma_s = -1 \Rightarrow  - \mathbf{E}_s - \mathbf{E}_s^T =- (\mathbf{E}_s + \mathbf{E}_s^T)  = \mathbf{0}$, therefore $ \mathbf{E}_s$ is anti-symmetric.
\end{proof}

As can be seen from the next examples the structure of the matrices is very symmetric.
The  diagonal holds the scalar value
\[
\pi (\grpart{a}_0) = \grpart{a}_0 \mathbf{I}
\]
while the co-diagonal holds the pseudoscalar value,

\section{Low-dimensional matrix representations}
\label{sec:examples}

In this section examples are given for non-degenerate Clifford algebras. 

\subsection{Representations for n=2}
\label{sec:n2}

For a general element of the form 
\[
A= {{a}_{1}}+{{e}_{1}}\,{{a}_{2}}+{{e}_{2}}\,{{a}_{3}}+{{a}_{4}}\,\left( {{e}_{1}}{{e}_{2}}\right) 
\]
we have the following canonical matrix representations:
\begin{itemize}
	\item Euclidean algebra  $\Cl{2,0}$	
	\[
	\mathbf{A}=\begin{pmatrix}{{a}_{1}} & {{a}_{2}} & {{a}_{3}} & {{a}_{4}}\\
	{{a}_{2}} & {{a}_{1}} & {{a}_{4}} & {{a}_{3}}\\
	{{a}_{3}} & -{{a}_{4}} & {{a}_{1}} & -{{a}_{2}}\\
	-{{a}_{4}} & {{a}_{3}} & -{{a}_{2}} & {{a}_{1}}
	\end{pmatrix}
	\]
	\item Split-quaternions $\Cl{1,1}$ \\
	\[
	\mathbf{A}=\begin{pmatrix}{{a}_{1}} & {{a}_{2}} & {{a}_{3}} & {{a}_{4}}\\
	{{a}_{2}} & {{a}_{1}} & {{a}_{4}} & {{a}_{3}}\\
	-{{a}_{3}} & {{a}_{4}} & {{a}_{1}} & -{{a}_{2}}\\
	{{a}_{4}} & -{{a}_{3}} & -{{a}_{2}} & {{a}_{1}}
	\end{pmatrix}
	\]
	\item Quaternions $\Cl{0,2}$ \\
	\[
	\mathbf{A}=\begin{pmatrix}{{a}_{1}} & {{a}_{2}} & {{a}_{3}} & {{a}_{4}}\\
	-{{a}_{2}} & {{a}_{1}} & -{{a}_{4}} & {{a}_{3}}\\
	-{{a}_{3}} & {{a}_{4}} & {{a}_{1}} & -{{a}_{2}}\\
	-{{a}_{4}} & -{{a}_{3}} & {{a}_{2}} & {{a}_{1}}
	\end{pmatrix}
	\]
\end{itemize}

\subsection{Representations for n=3}
\label{sec:n3}

For a general element of the form 
\[
{{a}_{1}}+{{e}_{1}}\,{{a}_{2}}+{{e}_{2}}\,{{a}_{3}}+{{e}_{3}}\,{{a}_{4}}+{{a}_{5}}\,\left( {{e}_{1}} {{e}_{2}}\right) 
+{{a}_{6}}\,\left( {{e}_{1}} {{e}_{3}}\right) +{{a}_{7}}\,\left( {{e}_{2}} {{e}_{3}}\right) +{{a}_{8}}\,\left( {{e}_{1}} {{e}_{2}} {{e}_{3}}\right)
\]
we have the following canonical matrix representations:

\begin{itemize}
	\item  Geometric algebra / Algebra of Physical Space $\Cl{3,0}$	
	\[
	\mathbf{A}=\begin{pmatrix}{{a}_{1}} & {{a}_{2}} & {{a}_{3}} & {{a}_{4}} & {{a}_{5}} & {{a}_{6}} & {{a}_{7}} & {{a}_{8}}\\
	{{a}_{2}} & {{a}_{1}} & {{a}_{5}} & {{a}_{6}} & {{a}_{3}} & {{a}_{4}} & {{a}_{8}} & {{a}_{7}}\\
	{{a}_{3}} & -{{a}_{5}} & {{a}_{1}} & {{a}_{7}} & -{{a}_{2}} & -{{a}_{8}} & {{a}_{4}} & -{{a}_{6}}\\
	{{a}_{4}} & -{{a}_{6}} & -{{a}_{7}} & {{a}_{1}} & {{a}_{8}} & -{{a}_{2}} & -{{a}_{3}} & {{a}_{5}}\\
	-{{a}_{5}} & {{a}_{3}} & -{{a}_{2}} & -{{a}_{8}} & {{a}_{1}} & {{a}_{7}} & -{{a}_{6}} & {{a}_{4}}\\
	-{{a}_{6}} & {{a}_{4}} & {{a}_{8}} & -{{a}_{2}} & -{{a}_{7}} & {{a}_{1}} & {{a}_{5}} & -{{a}_{3}}\\
	-{{a}_{7}} & -{{a}_{8}} & {{a}_{4}} & -{{a}_{3}} & {{a}_{6}} & -{{a}_{5}} & {{a}_{1}} & {{a}_{2}}\\
	-{{a}_{8}} & -{{a}_{7}} & {{a}_{6}} & -{{a}_{5}} & {{a}_{4}} & -{{a}_{3}} & {{a}_{2}} & {{a}_{1}}\end{pmatrix}
	\]
	
	\item Clifford algebra $\Cl{2,1}$	
	\[
	\mathbf{A}=\begin{pmatrix}{{a}_{1}} & {{a}_{2}} & {{a}_{3}} & {{a}_{4}} & {{a}_{5}} & {{a}_{6}} & {{a}_{7}} & {{a}_{8}}\\
	{{a}_{2}} & {{a}_{1}} & {{a}_{5}} & {{a}_{6}} & {{a}_{3}} & {{a}_{4}} & {{a}_{8}} & {{a}_{7}}\\
	{{a}_{3}} & -{{a}_{5}} & {{a}_{1}} & {{a}_{7}} & -{{a}_{2}} & -{{a}_{8}} & {{a}_{4}} & -{{a}_{6}}\\
	-{{a}_{4}} & {{a}_{6}} & {{a}_{7}} & {{a}_{1}} & -{{a}_{8}} & -{{a}_{2}} & -{{a}_{3}} & {{a}_{5}}\\
	-{{a}_{5}} & {{a}_{3}} & -{{a}_{2}} & -{{a}_{8}} & {{a}_{1}} & {{a}_{7}} & -{{a}_{6}} & {{a}_{4}}\\
	{{a}_{6}} & -{{a}_{4}} & -{{a}_{8}} & -{{a}_{2}} & {{a}_{7}} & {{a}_{1}} & {{a}_{5}} & -{{a}_{3}}\\
	{{a}_{7}} & {{a}_{8}} & -{{a}_{4}} & -{{a}_{3}} & -{{a}_{6}} & -{{a}_{5}} & {{a}_{1}} & {{a}_{2}}\\
	{{a}_{8}} & {{a}_{7}} & -{{a}_{6}} & -{{a}_{5}} & -{{a}_{4}} & -{{a}_{3}} & {{a}_{2}} & {{a}_{1}}\end{pmatrix}
	\]
	
	\item Clifford algebra $\Cl{1,2}$
	\[
	\mathbf{A}=\begin{pmatrix}{{a}_{1}} & {{a}_{2}} & {{a}_{3}} & {{a}_{4}} & {{a}_{5}} & {{a}_{6}} & {{a}_{7}} & {{a}_{8}}\\
	{{a}_{2}} & {{a}_{1}} & {{a}_{5}} & {{a}_{6}} & {{a}_{3}} & {{a}_{4}} & {{a}_{8}} & {{a}_{7}}\\
	-{{a}_{3}} & {{a}_{5}} & {{a}_{1}} & -{{a}_{7}} & -{{a}_{2}} & {{a}_{8}} & {{a}_{4}} & -{{a}_{6}}\\
	-{{a}_{4}} & {{a}_{6}} & {{a}_{7}} & {{a}_{1}} & -{{a}_{8}} & -{{a}_{2}} & -{{a}_{3}} & {{a}_{5}}\\
	{{a}_{5}} & -{{a}_{3}} & -{{a}_{2}} & {{a}_{8}} & {{a}_{1}} & -{{a}_{7}} & -{{a}_{6}} & {{a}_{4}}\\
	{{a}_{6}} & -{{a}_{4}} & -{{a}_{8}} & -{{a}_{2}} & {{a}_{7}} & {{a}_{1}} & {{a}_{5}} & -{{a}_{3}}\\
	-{{a}_{7}} & -{{a}_{8}} & -{{a}_{4}} & {{a}_{3}} & -{{a}_{6}} & {{a}_{5}} & {{a}_{1}} & {{a}_{2}}\\
	-{{a}_{8}} & -{{a}_{7}} & -{{a}_{6}} & {{a}_{5}} & -{{a}_{4}} & {{a}_{3}} & {{a}_{2}} & {{a}_{1}}\end{pmatrix}
	\]
	
	\item Clifford algebra $\Cl{0,3}$
	\[
	\mathbf{A}=\begin{pmatrix}{{a}_{1}} & {{a}_{2}} & {{a}_{3}} & {{a}_{4}} & {{a}_{5}} & {{a}_{6}} & {{a}_{7}} & {{a}_{8}}\\
	-{{a}_{2}} & {{a}_{1}} & -{{a}_{5}} & -{{a}_{6}} & {{a}_{3}} & {{a}_{4}} & -{{a}_{8}} & {{a}_{7}}\\
	-{{a}_{3}} & {{a}_{5}} & {{a}_{1}} & -{{a}_{7}} & -{{a}_{2}} & {{a}_{8}} & {{a}_{4}} & -{{a}_{6}}\\
	-{{a}_{4}} & {{a}_{6}} & {{a}_{7}} & {{a}_{1}} & -{{a}_{8}} & -{{a}_{2}} & -{{a}_{3}} & {{a}_{5}}\\
	-{{a}_{5}} & -{{a}_{3}} & {{a}_{2}} & -{{a}_{8}} & {{a}_{1}} & -{{a}_{7}} & {{a}_{6}} & {{a}_{4}}\\
	-{{a}_{6}} & -{{a}_{4}} & {{a}_{8}} & {{a}_{2}} & {{a}_{7}} & {{a}_{1}} & -{{a}_{5}} & -{{a}_{3}}\\
	-{{a}_{7}} & -{{a}_{8}} & -{{a}_{4}} & {{a}_{3}} & -{{a}_{6}} & {{a}_{5}} & {{a}_{1}} & {{a}_{2}}\\
	{{a}_{8}} & -{{a}_{7}} & {{a}_{6}} & -{{a}_{5}} & -{{a}_{4}} & {{a}_{3}} & -{{a}_{2}} & {{a}_{1}}\end{pmatrix}
	\]
	
\end{itemize}

\subsection{Representations for n=4}
\label{sec:n4}

For a general element of the form 
\begin{flalign*}
A= {{a}_{1}} +{{e}_{1}}\,{{a}_{2}}+{{e}_{2}}\,{{a}_{3}}+{{e}_{3}}\,{{a}_{4}}+{{e}_{4}}\,{{a}_{5}} +
{{a}_{6}}  \left( {{e}_{1}} {{e}_{2}}\right) +{{a}_{7}}\,\left( {{e}_{1}} {{e}_{3}}\right)  + \\
{{a}_{8}}\,\left( {{e}_{1}} {{e}_{4}}\right) +{{a}_{9}}\,\left( {{e}_{2}} {{e}_{3}}\right)  +{{a}_{10}}\,\left( {{e}_{2}} {{e}_{4}}\right) +{{a}_{11}}\,\left( {{e}_{3}} {{e}_{4}}\right) + \\
{{a}_{12}} \left( {{e}_{1}} {{e}_{2}} {{e}_{3}}\right) 
+{{a}_{13}}\,\left( {{e}_{1}} {{e}_{2}} {{e}_{4}}\right) +{{a}_{14}}\,\left( {{e}_{1}} {{e}_{3}} {{e}_{4}}\right) +{{a}_{15}}\,\left( {{e}_{2}} {{e}_{3}} {{e}_{4}}\right) +  {{a}_{16}}\,\left( {{e}_{1}} {{e}_{2}} {{e}_{3}} {{e}_{4}}\right)  
\end{flalign*}
we have the following canonical matrix representations:
\begin{itemize}
	\setcounter{MaxMatrixCols}{20}
	\setlength{\arraycolsep}{.5\arraycolsep}
	
	\item  Euclidean algebra $\Cl{4,0}$	
	\[\tiny
	\begin{pmatrix}{{a}_{1}} & {{a}_{2}} & {{a}_{3}} & {{a}_{4}} & {{a}_{5}} & {{a}_{6}} & {{a}_{7}} & {{a}_{8}} & {{a}_{9}} & {{a}_{10}} & {{a}_{11}} & {{a}_{12}} & {{a}_{13}} & {{a}_{14}} & {{a}_{15}} & {{a}_{16}}\\
	{{a}_{2}} & {{a}_{1}} & {{a}_{6}} & {{a}_{7}} & {{a}_{8}} & {{a}_{3}} & {{a}_{4}} & {{a}_{5}} & {{a}_{12}} & {{a}_{13}} & {{a}_{14}} & {{a}_{9}} & {{a}_{10}} & {{a}_{11}} & {{a}_{16}} & {{a}_{15}}\\
	{{a}_{3}} & -{{a}_{6}} & {{a}_{1}} & {{a}_{9}} & {{a}_{10}} & -{{a}_{2}} & -{{a}_{12}} & -{{a}_{13}} & {{a}_{4}} & {{a}_{5}} & {{a}_{15}} & -{{a}_{7}} & -{{a}_{8}} & -{{a}_{16}} & {{a}_{11}} & -{{a}_{14}}\\
	{{a}_{4}} & -{{a}_{7}} & -{{a}_{9}} & {{a}_{1}} & {{a}_{11}} & {{a}_{12}} & -{{a}_{2}} & -{{a}_{14}} & -{{a}_{3}} & -{{a}_{15}} & {{a}_{5}} & {{a}_{6}} & {{a}_{16}} & -{{a}_{8}} & -{{a}_{10}} & {{a}_{13}}\\
	{{a}_{5}} & -{{a}_{8}} & -{{a}_{10}} & -{{a}_{11}} & {{a}_{1}} & {{a}_{13}} & {{a}_{14}} & -{{a}_{2}} & {{a}_{15}} & -{{a}_{3}} & -{{a}_{4}} & -{{a}_{16}} & {{a}_{6}} & {{a}_{7}} & {{a}_{9}} & -{{a}_{12}}\\
	-{{a}_{6}} & {{a}_{3}} & -{{a}_{2}} & -{{a}_{12}} & -{{a}_{13}} & {{a}_{1}} & {{a}_{9}} & {{a}_{10}} & -{{a}_{7}} & -{{a}_{8}} & -{{a}_{16}} & {{a}_{4}} & {{a}_{5}} & {{a}_{15}} & -{{a}_{14}} & {{a}_{11}}\\
	-{{a}_{7}} & {{a}_{4}} & {{a}_{12}} & -{{a}_{2}} & -{{a}_{14}} & -{{a}_{9}} & {{a}_{1}} & {{a}_{11}} & {{a}_{6}} & {{a}_{16}} & -{{a}_{8}} & -{{a}_{3}} & -{{a}_{15}} & {{a}_{5}} & {{a}_{13}} & -{{a}_{10}}\\
	-{{a}_{8}} & {{a}_{5}} & {{a}_{13}} & {{a}_{14}} & -{{a}_{2}} & -{{a}_{10}} & -{{a}_{11}} & {{a}_{1}} & -{{a}_{16}} & {{a}_{6}} & {{a}_{7}} & {{a}_{15}} & -{{a}_{3}} & -{{a}_{4}} & -{{a}_{12}} & {{a}_{9}}\\
	-{{a}_{9}} & -{{a}_{12}} & {{a}_{4}} & -{{a}_{3}} & -{{a}_{15}} & {{a}_{7}} & -{{a}_{6}} & -{{a}_{16}} & {{a}_{1}} & {{a}_{11}} & -{{a}_{10}} & {{a}_{2}} & {{a}_{14}} & -{{a}_{13}} & {{a}_{5}} & {{a}_{8}}\\
	-{{a}_{10}} & -{{a}_{13}} & {{a}_{5}} & {{a}_{15}} & -{{a}_{3}} & {{a}_{8}} & {{a}_{16}} & -{{a}_{6}} & -{{a}_{11}} & {{a}_{1}} & {{a}_{9}} & -{{a}_{14}} & {{a}_{2}} & {{a}_{12}} & -{{a}_{4}} & -{{a}_{7}}\\
	-{{a}_{11}} & -{{a}_{14}} & -{{a}_{15}} & {{a}_{5}} & -{{a}_{4}} & -{{a}_{16}} & {{a}_{8}} & -{{a}_{7}} & {{a}_{10}} & -{{a}_{9}} & {{a}_{1}} & {{a}_{13}} & -{{a}_{12}} & {{a}_{2}} & {{a}_{3}} & {{a}_{6}}\\
	-{{a}_{12}} & -{{a}_{9}} & {{a}_{7}} & -{{a}_{6}} & -{{a}_{16}} & {{a}_{4}} & -{{a}_{3}} & -{{a}_{15}} & {{a}_{2}} & {{a}_{14}} & -{{a}_{13}} & {{a}_{1}} & {{a}_{11}} & -{{a}_{10}} & {{a}_{8}} & {{a}_{5}}\\
	-{{a}_{13}} & -{{a}_{10}} & {{a}_{8}} & {{a}_{16}} & -{{a}_{6}} & {{a}_{5}} & {{a}_{15}} & -{{a}_{3}} & -{{a}_{14}} & {{a}_{2}} & {{a}_{12}} & -{{a}_{11}} & {{a}_{1}} & {{a}_{9}} & -{{a}_{7}} & -{{a}_{4}}\\
	-{{a}_{14}} & -{{a}_{11}} & -{{a}_{16}} & {{a}_{8}} & -{{a}_{7}} & -{{a}_{15}} & {{a}_{5}} & -{{a}_{4}} & {{a}_{13}} & -{{a}_{12}} & {{a}_{2}} & {{a}_{10}} & -{{a}_{9}} & {{a}_{1}} & {{a}_{6}} & {{a}_{3}}\\
	-{{a}_{15}} & {{a}_{16}} & -{{a}_{11}} & {{a}_{10}} & -{{a}_{9}} & {{a}_{14}} & -{{a}_{13}} & {{a}_{12}} & {{a}_{5}} & -{{a}_{4}} & {{a}_{3}} & -{{a}_{8}} & {{a}_{7}} & -{{a}_{6}} & {{a}_{1}} & -{{a}_{2}}\\
	{{a}_{16}} & -{{a}_{15}} & {{a}_{14}} & -{{a}_{13}} & {{a}_{12}} & -{{a}_{11}} & {{a}_{10}} & -{{a}_{9}} & -{{a}_{8}} & {{a}_{7}} & -{{a}_{6}} & {{a}_{5}} & -{{a}_{4}} & {{a}_{3}} & -{{a}_{2}} & {{a}_{1}}\end{pmatrix}
	\]
	
	\item Clifford algebra $\Cl{3,1}$
	\[\tiny
	\begin{pmatrix}
	{{a}_{1}} & {{a}_{2}} & {{a}_{3}} & {{a}_{4}} & {{a}_{5}} & {{a}_{6}} & {{a}_{7}} & {{a}_{8}} & {{a}_{9}} & {{a}_{10}} & {{a}_{11}} & {{a}_{12}} & {{a}_{13}} & {{a}_{14}} & {{a}_{15}} & {{a}_{16}}\\
	{{a}_{2}} & {{a}_{1}} & {{a}_{6}} & {{a}_{7}} & {{a}_{8}} & {{a}_{3}} & {{a}_{4}} & {{a}_{5}} & {{a}_{12}} & {{a}_{13}} & {{a}_{14}} & {{a}_{9}} & {{a}_{10}} & {{a}_{11}} & {{a}_{16}} & {{a}_{15}}\\
	{{a}_{3}} & -{{a}_{6}} & {{a}_{1}} & {{a}_{9}} & {{a}_{10}} & -{{a}_{2}} & -{{a}_{12}} & -{{a}_{13}} & {{a}_{4}} & {{a}_{5}} & {{a}_{15}} & -{{a}_{7}} & -{{a}_{8}} & -{{a}_{16}} & {{a}_{11}} & -{{a}_{14}}\\
	{{a}_{4}} & -{{a}_{7}} & -{{a}_{9}} & {{a}_{1}} & {{a}_{11}} & {{a}_{12}} & -{{a}_{2}} & -{{a}_{14}} & -{{a}_{3}} & -{{a}_{15}} & {{a}_{5}} & {{a}_{6}} & {{a}_{16}} & -{{a}_{8}} & -{{a}_{10}} & {{a}_{13}}\\
	{{a}_{5}} & -{{a}_{8}} & -{{a}_{10}} & -{{a}_{11}} & {{a}_{1}} & {{a}_{13}} & {{a}_{14}} & -{{a}_{2}} & {{a}_{15}} & -{{a}_{3}} & -{{a}_{4}} & -{{a}_{16}} & {{a}_{6}} & {{a}_{7}} & {{a}_{9}} & -{{a}_{12}}\\
	-{{a}_{6}} & {{a}_{3}} & -{{a}_{2}} & -{{a}_{12}} & -{{a}_{13}} & {{a}_{1}} & {{a}_{9}} & {{a}_{10}} & -{{a}_{7}} & -{{a}_{8}} & -{{a}_{16}} & {{a}_{4}} & {{a}_{5}} & {{a}_{15}} & -{{a}_{14}} & {{a}_{11}}\\
	-{{a}_{7}} & {{a}_{4}} & {{a}_{12}} & -{{a}_{2}} & -{{a}_{14}} & -{{a}_{9}} & {{a}_{1}} & {{a}_{11}} & {{a}_{6}} & {{a}_{16}} & -{{a}_{8}} & -{{a}_{3}} & -{{a}_{15}} & {{a}_{5}} & {{a}_{13}} & -{{a}_{10}}\\
	-{{a}_{8}} & {{a}_{5}} & {{a}_{13}} & {{a}_{14}} & -{{a}_{2}} & -{{a}_{10}} & -{{a}_{11}} & {{a}_{1}} & -{{a}_{16}} & {{a}_{6}} & {{a}_{7}} & {{a}_{15}} & -{{a}_{3}} & -{{a}_{4}} & -{{a}_{12}} & {{a}_{9}}\\
	-{{a}_{9}} & -{{a}_{12}} & {{a}_{4}} & -{{a}_{3}} & -{{a}_{15}} & {{a}_{7}} & -{{a}_{6}} & -{{a}_{16}} & {{a}_{1}} & {{a}_{11}} & -{{a}_{10}} & {{a}_{2}} & {{a}_{14}} & -{{a}_{13}} & {{a}_{5}} & {{a}_{8}}\\
	-{{a}_{10}} & -{{a}_{13}} & {{a}_{5}} & {{a}_{15}} & -{{a}_{3}} & {{a}_{8}} & {{a}_{16}} & -{{a}_{6}} & -{{a}_{11}} & {{a}_{1}} & {{a}_{9}} & -{{a}_{14}} & {{a}_{2}} & {{a}_{12}} & -{{a}_{4}} & -{{a}_{7}}\\
	-{{a}_{11}} & -{{a}_{14}} & -{{a}_{15}} & {{a}_{5}} & -{{a}_{4}} & -{{a}_{16}} & {{a}_{8}} & -{{a}_{7}} & {{a}_{10}} & -{{a}_{9}} & {{a}_{1}} & {{a}_{13}} & -{{a}_{12}} & {{a}_{2}} & {{a}_{3}} & {{a}_{6}}\\
	-{{a}_{12}} & -{{a}_{9}} & {{a}_{7}} & -{{a}_{6}} & -{{a}_{16}} & {{a}_{4}} & -{{a}_{3}} & -{{a}_{15}} & {{a}_{2}} & {{a}_{14}} & -{{a}_{13}} & {{a}_{1}} & {{a}_{11}} & -{{a}_{10}} & {{a}_{8}} & {{a}_{5}}\\
	-{{a}_{13}} & -{{a}_{10}} & {{a}_{8}} & {{a}_{16}} & -{{a}_{6}} & {{a}_{5}} & {{a}_{15}} & -{{a}_{3}} & -{{a}_{14}} & {{a}_{2}} & {{a}_{12}} & -{{a}_{11}} & {{a}_{1}} & {{a}_{9}} & -{{a}_{7}} & -{{a}_{4}}\\
	-{{a}_{14}} & -{{a}_{11}} & -{{a}_{16}} & {{a}_{8}} & -{{a}_{7}} & -{{a}_{15}} & {{a}_{5}} & -{{a}_{4}} & {{a}_{13}} & -{{a}_{12}} & {{a}_{2}} & {{a}_{10}} & -{{a}_{9}} & {{a}_{1}} & {{a}_{6}} & {{a}_{3}}\\
	-{{a}_{15}} & {{a}_{16}} & -{{a}_{11}} & {{a}_{10}} & -{{a}_{9}} & {{a}_{14}} & -{{a}_{13}} & {{a}_{12}} & {{a}_{5}} & -{{a}_{4}} & {{a}_{3}} & -{{a}_{8}} & {{a}_{7}} & -{{a}_{6}} & {{a}_{1}} & -{{a}_{2}}\\
	{{a}_{16}} & -{{a}_{15}} & {{a}_{14}} & -{{a}_{13}} & {{a}_{12}} & -{{a}_{11}} & {{a}_{10}} & -{{a}_{9}} & -{{a}_{8}} & {{a}_{7}} & -{{a}_{6}} & {{a}_{5}} & -{{a}_{4}} & {{a}_{3}} & -{{a}_{2}} & {{a}_{1}}
	\end{pmatrix}
	\]
	
	\item Clifford algebra $\Cl{2,2}$	
	\[\tiny
	\begin{pmatrix}{{a}_{1}} & {{a}_{2}} & {{a}_{3}} & {{a}_{4}} & {{a}_{5}} & {{a}_{6}} & {{a}_{7}} & {{a}_{8}} & {{a}_{9}} & {{a}_{10}} & {{a}_{11}} & {{a}_{12}} & {{a}_{13}} & {{a}_{14}} & {{a}_{15}} & {{a}_{16}}\\
	{{a}_{2}} & {{a}_{1}} & {{a}_{6}} & {{a}_{7}} & {{a}_{8}} & {{a}_{3}} & {{a}_{4}} & {{a}_{5}} & {{a}_{12}} & {{a}_{13}} & {{a}_{14}} & {{a}_{9}} & {{a}_{10}} & {{a}_{11}} & {{a}_{16}} & {{a}_{15}}\\
	{{a}_{3}} & -{{a}_{6}} & {{a}_{1}} & {{a}_{9}} & {{a}_{10}} & -{{a}_{2}} & -{{a}_{12}} & -{{a}_{13}} & {{a}_{4}} & {{a}_{5}} & {{a}_{15}} & -{{a}_{7}} & -{{a}_{8}} & -{{a}_{16}} & {{a}_{11}} & -{{a}_{14}}\\
	-{{a}_{4}} & {{a}_{7}} & {{a}_{9}} & {{a}_{1}} & -{{a}_{11}} & -{{a}_{12}} & -{{a}_{2}} & {{a}_{14}} & -{{a}_{3}} & {{a}_{15}} & {{a}_{5}} & {{a}_{6}} & -{{a}_{16}} & -{{a}_{8}} & -{{a}_{10}} & {{a}_{13}}\\
	-{{a}_{5}} & {{a}_{8}} & {{a}_{10}} & {{a}_{11}} & {{a}_{1}} & -{{a}_{13}} & -{{a}_{14}} & -{{a}_{2}} & -{{a}_{15}} & -{{a}_{3}} & -{{a}_{4}} & {{a}_{16}} & {{a}_{6}} & {{a}_{7}} & {{a}_{9}} & -{{a}_{12}}\\
	-{{a}_{6}} & {{a}_{3}} & -{{a}_{2}} & -{{a}_{12}} & -{{a}_{13}} & {{a}_{1}} & {{a}_{9}} & {{a}_{10}} & -{{a}_{7}} & -{{a}_{8}} & -{{a}_{16}} & {{a}_{4}} & {{a}_{5}} & {{a}_{15}} & -{{a}_{14}} & {{a}_{11}}\\
	{{a}_{7}} & -{{a}_{4}} & -{{a}_{12}} & -{{a}_{2}} & {{a}_{14}} & {{a}_{9}} & {{a}_{1}} & -{{a}_{11}} & {{a}_{6}} & -{{a}_{16}} & -{{a}_{8}} & -{{a}_{3}} & {{a}_{15}} & {{a}_{5}} & {{a}_{13}} & -{{a}_{10}}\\
	{{a}_{8}} & -{{a}_{5}} & -{{a}_{13}} & -{{a}_{14}} & -{{a}_{2}} & {{a}_{10}} & {{a}_{11}} & {{a}_{1}} & {{a}_{16}} & {{a}_{6}} & {{a}_{7}} & -{{a}_{15}} & -{{a}_{3}} & -{{a}_{4}} & -{{a}_{12}} & {{a}_{9}}\\
	{{a}_{9}} & {{a}_{12}} & -{{a}_{4}} & -{{a}_{3}} & {{a}_{15}} & -{{a}_{7}} & -{{a}_{6}} & {{a}_{16}} & {{a}_{1}} & -{{a}_{11}} & -{{a}_{10}} & {{a}_{2}} & -{{a}_{14}} & -{{a}_{13}} & {{a}_{5}} & {{a}_{8}}\\
	{{a}_{10}} & {{a}_{13}} & -{{a}_{5}} & -{{a}_{15}} & -{{a}_{3}} & -{{a}_{8}} & -{{a}_{16}} & -{{a}_{6}} & {{a}_{11}} & {{a}_{1}} & {{a}_{9}} & {{a}_{14}} & {{a}_{2}} & {{a}_{12}} & -{{a}_{4}} & -{{a}_{7}}\\
	-{{a}_{11}} & -{{a}_{14}} & -{{a}_{15}} & -{{a}_{5}} & {{a}_{4}} & -{{a}_{16}} & -{{a}_{8}} & {{a}_{7}} & -{{a}_{10}} & {{a}_{9}} & {{a}_{1}} & -{{a}_{13}} & {{a}_{12}} & {{a}_{2}} & {{a}_{3}} & {{a}_{6}}\\
	{{a}_{12}} & {{a}_{9}} & -{{a}_{7}} & -{{a}_{6}} & {{a}_{16}} & -{{a}_{4}} & -{{a}_{3}} & {{a}_{15}} & {{a}_{2}} & -{{a}_{14}} & -{{a}_{13}} & {{a}_{1}} & -{{a}_{11}} & -{{a}_{10}} & {{a}_{8}} & {{a}_{5}}\\
	{{a}_{13}} & {{a}_{10}} & -{{a}_{8}} & -{{a}_{16}} & -{{a}_{6}} & -{{a}_{5}} & -{{a}_{15}} & -{{a}_{3}} & {{a}_{14}} & {{a}_{2}} & {{a}_{12}} & {{a}_{11}} & {{a}_{1}} & {{a}_{9}} & -{{a}_{7}} & -{{a}_{4}}\\
	-{{a}_{14}} & -{{a}_{11}} & -{{a}_{16}} & -{{a}_{8}} & {{a}_{7}} & -{{a}_{15}} & -{{a}_{5}} & {{a}_{4}} & -{{a}_{13}} & {{a}_{12}} & {{a}_{2}} & -{{a}_{10}} & {{a}_{9}} & {{a}_{1}} & {{a}_{6}} & {{a}_{3}}\\
	-{{a}_{15}} & {{a}_{16}} & -{{a}_{11}} & -{{a}_{10}} & {{a}_{9}} & {{a}_{14}} & {{a}_{13}} & -{{a}_{12}} & -{{a}_{5}} & {{a}_{4}} & {{a}_{3}} & {{a}_{8}} & -{{a}_{7}} & -{{a}_{6}} & {{a}_{1}} & -{{a}_{2}}\\
	{{a}_{16}} & -{{a}_{15}} & {{a}_{14}} & {{a}_{13}} & -{{a}_{12}} & -{{a}_{11}} & -{{a}_{10}} & {{a}_{9}} & {{a}_{8}} & -{{a}_{7}} & -{{a}_{6}} & -{{a}_{5}} & {{a}_{4}} & {{a}_{3}} & -{{a}_{2}} & {{a}_{1}}\end{pmatrix}
	\]
	
	\item  Space-time algebra $\Cl{1,3}$
	\[\tiny
	\begin{pmatrix}{{a}_{1}} & {{a}_{2}} & {{a}_{3}} & {{a}_{4}} & {{a}_{5}} & {{a}_{6}} & {{a}_{7}} & {{a}_{8}} & {{a}_{9}} & {{a}_{10}} & {{a}_{11}} & {{a}_{12}} & {{a}_{13}} & {{a}_{14}} & {{a}_{15}} & {{a}_{16}}\\
	{{a}_{2}} & {{a}_{1}} & {{a}_{6}} & {{a}_{7}} & {{a}_{8}} & {{a}_{3}} & {{a}_{4}} & {{a}_{5}} & {{a}_{12}} & {{a}_{13}} & {{a}_{14}} & {{a}_{9}} & {{a}_{10}} & {{a}_{11}} & {{a}_{16}} & {{a}_{15}}\\
	-{{a}_{3}} & {{a}_{6}} & {{a}_{1}} & -{{a}_{9}} & -{{a}_{10}} & -{{a}_{2}} & {{a}_{12}} & {{a}_{13}} & {{a}_{4}} & {{a}_{5}} & -{{a}_{15}} & -{{a}_{7}} & -{{a}_{8}} & {{a}_{16}} & {{a}_{11}} & -{{a}_{14}}\\
	-{{a}_{4}} & {{a}_{7}} & {{a}_{9}} & {{a}_{1}} & -{{a}_{11}} & -{{a}_{12}} & -{{a}_{2}} & {{a}_{14}} & -{{a}_{3}} & {{a}_{15}} & {{a}_{5}} & {{a}_{6}} & -{{a}_{16}} & -{{a}_{8}} & -{{a}_{10}} & {{a}_{13}}\\
	-{{a}_{5}} & {{a}_{8}} & {{a}_{10}} & {{a}_{11}} & {{a}_{1}} & -{{a}_{13}} & -{{a}_{14}} & -{{a}_{2}} & -{{a}_{15}} & -{{a}_{3}} & -{{a}_{4}} & {{a}_{16}} & {{a}_{6}} & {{a}_{7}} & {{a}_{9}} & -{{a}_{12}}\\
	{{a}_{6}} & -{{a}_{3}} & -{{a}_{2}} & {{a}_{12}} & {{a}_{13}} & {{a}_{1}} & -{{a}_{9}} & -{{a}_{10}} & -{{a}_{7}} & -{{a}_{8}} & {{a}_{16}} & {{a}_{4}} & {{a}_{5}} & -{{a}_{15}} & -{{a}_{14}} & {{a}_{11}}\\
	{{a}_{7}} & -{{a}_{4}} & -{{a}_{12}} & -{{a}_{2}} & {{a}_{14}} & {{a}_{9}} & {{a}_{1}} & -{{a}_{11}} & {{a}_{6}} & -{{a}_{16}} & -{{a}_{8}} & -{{a}_{3}} & {{a}_{15}} & {{a}_{5}} & {{a}_{13}} & -{{a}_{10}}\\
	{{a}_{8}} & -{{a}_{5}} & -{{a}_{13}} & -{{a}_{14}} & -{{a}_{2}} & {{a}_{10}} & {{a}_{11}} & {{a}_{1}} & {{a}_{16}} & {{a}_{6}} & {{a}_{7}} & -{{a}_{15}} & -{{a}_{3}} & -{{a}_{4}} & -{{a}_{12}} & {{a}_{9}}\\
	-{{a}_{9}} & -{{a}_{12}} & -{{a}_{4}} & {{a}_{3}} & -{{a}_{15}} & -{{a}_{7}} & {{a}_{6}} & -{{a}_{16}} & {{a}_{1}} & -{{a}_{11}} & {{a}_{10}} & {{a}_{2}} & -{{a}_{14}} & {{a}_{13}} & {{a}_{5}} & {{a}_{8}}\\
	-{{a}_{10}} & -{{a}_{13}} & -{{a}_{5}} & {{a}_{15}} & {{a}_{3}} & -{{a}_{8}} & {{a}_{16}} & {{a}_{6}} & {{a}_{11}} & {{a}_{1}} & -{{a}_{9}} & {{a}_{14}} & {{a}_{2}} & -{{a}_{12}} & -{{a}_{4}} & -{{a}_{7}}\\
	-{{a}_{11}} & -{{a}_{14}} & -{{a}_{15}} & -{{a}_{5}} & {{a}_{4}} & -{{a}_{16}} & -{{a}_{8}} & {{a}_{7}} & -{{a}_{10}} & {{a}_{9}} & {{a}_{1}} & -{{a}_{13}} & {{a}_{12}} & {{a}_{2}} & {{a}_{3}} & {{a}_{6}}\\
	-{{a}_{12}} & -{{a}_{9}} & -{{a}_{7}} & {{a}_{6}} & -{{a}_{16}} & -{{a}_{4}} & {{a}_{3}} & -{{a}_{15}} & {{a}_{2}} & -{{a}_{14}} & {{a}_{13}} & {{a}_{1}} & -{{a}_{11}} & {{a}_{10}} & {{a}_{8}} & {{a}_{5}}\\
	-{{a}_{13}} & -{{a}_{10}} & -{{a}_{8}} & {{a}_{16}} & {{a}_{6}} & -{{a}_{5}} & {{a}_{15}} & {{a}_{3}} & {{a}_{14}} & {{a}_{2}} & -{{a}_{12}} & {{a}_{11}} & {{a}_{1}} & -{{a}_{9}} & -{{a}_{7}} & -{{a}_{4}}\\
	-{{a}_{14}} & -{{a}_{11}} & -{{a}_{16}} & -{{a}_{8}} & {{a}_{7}} & -{{a}_{15}} & -{{a}_{5}} & {{a}_{4}} & -{{a}_{13}} & {{a}_{12}} & {{a}_{2}} & -{{a}_{10}} & {{a}_{9}} & {{a}_{1}} & {{a}_{6}} & {{a}_{3}}\\
	{{a}_{15}} & -{{a}_{16}} & -{{a}_{11}} & {{a}_{10}} & -{{a}_{9}} & {{a}_{14}} & -{{a}_{13}} & {{a}_{12}} & -{{a}_{5}} & {{a}_{4}} & -{{a}_{3}} & {{a}_{8}} & -{{a}_{7}} & {{a}_{6}} & {{a}_{1}} & -{{a}_{2}}\\
	-{{a}_{16}} & {{a}_{15}} & {{a}_{14}} & -{{a}_{13}} & {{a}_{12}} & -{{a}_{11}} & {{a}_{10}} & -{{a}_{9}} & {{a}_{8}} & -{{a}_{7}} & {{a}_{6}} & -{{a}_{5}} & {{a}_{4}} & -{{a}_{3}} & -{{a}_{2}} & {{a}_{1}}\end{pmatrix}
	\]
	
	\item Clifford algebra $\Cl{0,4}$
	\[\tiny
	\begin{pmatrix}{{a}_{1}} & {{a}_{2}} & {{a}_{3}} & {{a}_{4}} & {{a}_{5}} & {{a}_{6}} & {{a}_{7}} & {{a}_{8}} & {{a}_{9}} & {{a}_{10}} & {{a}_{11}} & {{a}_{12}} & {{a}_{13}} & {{a}_{14}} & {{a}_{15}} & {{a}_{16}}\\
	-{{a}_{2}} & {{a}_{1}} & -{{a}_{6}} & -{{a}_{7}} & -{{a}_{8}} & {{a}_{3}} & {{a}_{4}} & {{a}_{5}} & -{{a}_{12}} & -{{a}_{13}} & -{{a}_{14}} & {{a}_{9}} & {{a}_{10}} & {{a}_{11}} & -{{a}_{16}} & {{a}_{15}}\\
	-{{a}_{3}} & {{a}_{6}} & {{a}_{1}} & -{{a}_{9}} & -{{a}_{10}} & -{{a}_{2}} & {{a}_{12}} & {{a}_{13}} & {{a}_{4}} & {{a}_{5}} & -{{a}_{15}} & -{{a}_{7}} & -{{a}_{8}} & {{a}_{16}} & {{a}_{11}} & -{{a}_{14}}\\
	-{{a}_{4}} & {{a}_{7}} & {{a}_{9}} & {{a}_{1}} & -{{a}_{11}} & -{{a}_{12}} & -{{a}_{2}} & {{a}_{14}} & -{{a}_{3}} & {{a}_{15}} & {{a}_{5}} & {{a}_{6}} & -{{a}_{16}} & -{{a}_{8}} & -{{a}_{10}} & {{a}_{13}}\\
	-{{a}_{5}} & {{a}_{8}} & {{a}_{10}} & {{a}_{11}} & {{a}_{1}} & -{{a}_{13}} & -{{a}_{14}} & -{{a}_{2}} & -{{a}_{15}} & -{{a}_{3}} & -{{a}_{4}} & {{a}_{16}} & {{a}_{6}} & {{a}_{7}} & {{a}_{9}} & -{{a}_{12}}\\
	-{{a}_{6}} & -{{a}_{3}} & {{a}_{2}} & -{{a}_{12}} & -{{a}_{13}} & {{a}_{1}} & -{{a}_{9}} & -{{a}_{10}} & {{a}_{7}} & {{a}_{8}} & -{{a}_{16}} & {{a}_{4}} & {{a}_{5}} & -{{a}_{15}} & {{a}_{14}} & {{a}_{11}}\\
	-{{a}_{7}} & -{{a}_{4}} & {{a}_{12}} & {{a}_{2}} & -{{a}_{14}} & {{a}_{9}} & {{a}_{1}} & -{{a}_{11}} & -{{a}_{6}} & {{a}_{16}} & {{a}_{8}} & -{{a}_{3}} & {{a}_{15}} & {{a}_{5}} & -{{a}_{13}} & -{{a}_{10}}\\
	-{{a}_{8}} & -{{a}_{5}} & {{a}_{13}} & {{a}_{14}} & {{a}_{2}} & {{a}_{10}} & {{a}_{11}} & {{a}_{1}} & -{{a}_{16}} & -{{a}_{6}} & -{{a}_{7}} & -{{a}_{15}} & -{{a}_{3}} & -{{a}_{4}} & {{a}_{12}} & {{a}_{9}}\\
	-{{a}_{9}} & -{{a}_{12}} & -{{a}_{4}} & {{a}_{3}} & -{{a}_{15}} & -{{a}_{7}} & {{a}_{6}} & -{{a}_{16}} & {{a}_{1}} & -{{a}_{11}} & {{a}_{10}} & {{a}_{2}} & -{{a}_{14}} & {{a}_{13}} & {{a}_{5}} & {{a}_{8}}\\
	-{{a}_{10}} & -{{a}_{13}} & -{{a}_{5}} & {{a}_{15}} & {{a}_{3}} & -{{a}_{8}} & {{a}_{16}} & {{a}_{6}} & {{a}_{11}} & {{a}_{1}} & -{{a}_{9}} & {{a}_{14}} & {{a}_{2}} & -{{a}_{12}} & -{{a}_{4}} & -{{a}_{7}}\\
	-{{a}_{11}} & -{{a}_{14}} & -{{a}_{15}} & -{{a}_{5}} & {{a}_{4}} & -{{a}_{16}} & -{{a}_{8}} & {{a}_{7}} & -{{a}_{10}} & {{a}_{9}} & {{a}_{1}} & -{{a}_{13}} & {{a}_{12}} & {{a}_{2}} & {{a}_{3}} & {{a}_{6}}\\
	{{a}_{12}} & -{{a}_{9}} & {{a}_{7}} & -{{a}_{6}} & {{a}_{16}} & -{{a}_{4}} & {{a}_{3}} & -{{a}_{15}} & -{{a}_{2}} & {{a}_{14}} & -{{a}_{13}} & {{a}_{1}} & -{{a}_{11}} & {{a}_{10}} & -{{a}_{8}} & {{a}_{5}}\\
	{{a}_{13}} & -{{a}_{10}} & {{a}_{8}} & -{{a}_{16}} & -{{a}_{6}} & -{{a}_{5}} & {{a}_{15}} & {{a}_{3}} & -{{a}_{14}} & -{{a}_{2}} & {{a}_{12}} & {{a}_{11}} & {{a}_{1}} & -{{a}_{9}} & {{a}_{7}} & -{{a}_{4}}\\
	{{a}_{14}} & -{{a}_{11}} & {{a}_{16}} & {{a}_{8}} & -{{a}_{7}} & -{{a}_{15}} & -{{a}_{5}} & {{a}_{4}} & {{a}_{13}} & -{{a}_{12}} & -{{a}_{2}} & -{{a}_{10}} & {{a}_{9}} & {{a}_{1}} & -{{a}_{6}} & {{a}_{3}}\\
	{{a}_{15}} & -{{a}_{16}} & -{{a}_{11}} & {{a}_{10}} & -{{a}_{9}} & {{a}_{14}} & -{{a}_{13}} & {{a}_{12}} & -{{a}_{5}} & {{a}_{4}} & -{{a}_{3}} & {{a}_{8}} & -{{a}_{7}} & {{a}_{6}} & {{a}_{1}} & -{{a}_{2}}\\
	{{a}_{16}} & {{a}_{15}} & -{{a}_{14}} & {{a}_{13}} & -{{a}_{12}} & -{{a}_{11}} & {{a}_{10}} & -{{a}_{9}} & -{{a}_{8}} & {{a}_{7}} & -{{a}_{6}} & -{{a}_{5}} & {{a}_{4}} & -{{a}_{3}} & {{a}_{2}} & {{a}_{1}}\end{pmatrix}
	\]
	
\end{itemize}

Computations have been carried out using the computer code included in the Appendix.

\section{Complex matrix representation of central algebras}\label{sec:comprep}

In the present section we exhibit an algorithm for the construction of the complex matrix representation of central Clifford algebras. 
As an example, consider the dual representation of the \Cl{3,0}, where $I^2=-1$.
The matrix multiplication table can be represented as
\[\small
\mathbf{M}= \left( \begin{array}{cccc|cccc}
	1 & {e_1} & {e_2} & {e_3} & i {e_3} & -i {e_2} & i {e_1} & i\\
	{e_1} & 1 & {e_1}{e_2} & {e_1}{e_3} & i \left( {e_1}{e_3}\right)  & -i \left( {e_1}{e_2}\right)  & i & i {e_1}\\
	{e_2} & -{e_1}{e_2} & 1 & {e_2}{e_3} & i \left( {e_2}{e_3}\right)  & -i & -i \left( {e_1}{e_2}\right)  & i {e_2}\\
	{e_3} & -{e_1}{e_3} & -{e_2}{e_3} & 1 & i & i \left( {e_2}{e_3}\right)  & -i \left( {e_1}{e_3}\right)  & i {e_3}\\
	\hline
	i {e_3} & -i \left( {e_1}{e_3}\right)  & -i \left( {e_2}{e_3}\right)  & i & -1 & -{e_2}{e_3} & {e_1}{e_3} & -{e_3}\\
	-i {e_2} & i \left( {e_1}{e_2}\right)  & -i & -i \left( {e_2}{e_3}\right)  & {e_2}{e_3} & -1 & -{e_1}{e_2} & {e_2}\\
	i {e_1} & i & i \left( {e_1}{e_2}\right)  & i \left( {e_1}{e_3}\right)  & -{e_1}{e_3} & {e_1}{e_2} & -1 & -{e_1}\\
	i & i {e_1} & i {e_2} & i {e_3} & -{e_3} & {e_2} & -{e_1} & -1
\end{array}\right) 
\]
where we have substituted the pseudoscalar $I$ for the imaginary scalar unit $i$.
The main diagonal blocks do not contain \textit{i}, while the co-diagoanl ones do.
$\mathbf{M}$ can be further transformed into
\[
\mathbf{M}^{*}=
\begin{pmatrix}
	\mathbf{M}_{11} & i \mathbf{M}_{12} \\
	i \mathbf{M}_{21} & - \mathbf{M}_{22} 
\end{pmatrix}
=
\left( \begin{array}{cccc|cccc}
	1 & {e_1} & {e_2} & {e_3} & i {e_3} & -i {e_2} & i {e_1} & i\\
	{e_1} & 1 & i {e_3} & -i {e_2} & {e_2} & {e_3} & i & i {e_1}\\
	{e_2} & -i {e_3} & 1 & i {e_1} & -{e_1} & -i & {e_3} & i {e_2}\\
	{e_3} & i {e_2} & -i {e_1} & 1 & i & -{e_1} & -{e_2} & i {e_3}\\
	\hline
	i {e_3} & -{e_2} & {e_1} & i & -1 & -i {e_1} & -i {e_2} & -{e_3}\\
	-i {e_2} & -{e_3} & -i & {e_1} & i {e_1} & -1 & -i {e_3} & {e_2}\\
	i {e_1} & i & -{e_3} & {e_2} & i {e_2} & i {e_3} & -1 & -{e_1}\\
	i & i {e_1} & i {e_2} & i {e_3} & -{e_3} & {e_2} & -{e_1} & -1
\end{array}\right) 
\]
where this structure is not obvious, but we are left only with the generators.
Therefore, the same representation algorithm can be applied.
We use the complexified set of generators
\[
iE= \{  i {e_3}, -i {e_2}, i {e_1} , i\}
\]
and compute the reduced multiplication tables $-\mathbf{M}_{22} = \mathbf{Mat} (iE)$ and $\mathbf{G}= diag \left\langle - \mathbf{M}_{22}\right\rangle $.
The general algebra element can be represented as
\[
A= a +v + I (b + p)
\] 
Therefore, by linearity $\pi (A) =  \pi(a+ v) + \pi(I) \pi(b + p) $.
Therefore, we can make a substitution $ \pi(I) \mapsto i$ and proceed as before.
In such way, the complex matrix representation will be given by 
\[
\mathbf{ A}=\sum_{s \in \operatorname{ord} E} a_s [ \operatorname{diag} \mathbf{M}_{22}] _{s}C_s ( \mathbf{M}_{22} ) 
\]
Therefore, one can arrive at the following representation:
\[
\mathbf{ A}= \begin{pmatrix}{a_1}-i {a_8} & i {a_2}+{a_7} & i {a_3}-{a_6} & {a_4}-i {a_5}\\
	-i {a_2}-{a_7} & {a_1}-i {a_8} & i {a_4}+{a_5} & -{a_3}-i {a_6}\\
	-i {a_3}+{a_6} & -i {a_4}-{a_5} & {a_1}-i {a_8} & {a_2}-i {a_7}\\
	{a_4}-i {a_5} & -{a_3}-i {a_6} & {a_2}-i {a_7} & {a_1}-i {a_8}
\end{pmatrix}\]
The advantage of this representation is that the trace produces the pseudoscalar plus the scalar components in an obvious manner
$
\mathrm{Tr} \mathbf{ A} =  4 \left({a_1}-i {a_8} \right) 
$.

\section{Computation of multivector inverses}\label{sec:appli}

Computation of Clifford inverse has drawn attention in the literature \cite{Acus2018, Hitzer2019, Shirokov2020}.
Multivector inverses can be computed using the matrix representation and the characteristic polynomial.
The matrix inverse is 
\[
\mathbf{A}^{-1} = \frac{1}{|| \mathbf{A} ||} \mathbf{A}^{\#}
\]
where $\#$ denotes the matrix adjunct operation and $|| \mathbf{A} ||$ denotes the matrix determinant.
The formula is not practical, because it requires the computation of  $n^2+1$ determinants. 
With the help of the Cayley-Hamilton Theorem, the inverse of \textbf{A} can be expressed as a polynomial in \textbf{A}.
The inverse can be computed as the last step of the  Faddeev–LeVerrier--Souriau (FVS) algorithm \cite{Souriau1948, Faddeev1949}.
The algorithm computes the coefficients of the characteristic  polynomial
\[
p(\lambda):= ||\lambda I_n  - \mathbf{A} || = \lambda^n + c_1 \lambda^{n-1} + \ldots c_{n-1} \lambda + c_n, \quad n= dim(\mathbf{A})
\]
in \textit{n} steps:
\[
\begin{array}{l|lr}
	\mathbf{M}_1 = \mathbf{A}                                   & t_1 = \ \  \mathrm{Tr} [\mathbf{M}_1]            & c_1=- t_1 \\
	\mathbf{M}_2 = \mathbf{A} \mathbf{M}_1 - t_1 \mathbf{I}     & t_2 =\frac{1}{2}\mathrm{Tr}[\mathbf{A} \mathbf{M}_1] & c_2=- t_2 \\
	$\ldots$ & $\ldots$ & $\ldots$ \\
	\mathbf{M}_k = \mathbf{A} \mathbf{M}_{k-1} - t_k \mathbf{I} & t_k =\frac{1}{k}\mathrm{Tr}[\mathbf{A} \mathbf{M}_{k-1}] & c_k=- t_k
\end{array}
\]
with $\mathbf{M}_n =\mathbf{0} $.
The matrix inverse can be computed from the last step of the algorithm as
\[
\mathbf{A}^{-1} = \frac{1}{t_n} \mathbf{M}_{n-1}
\]
under the obvious restriction $t_n \neq 0$.

\begin{example}
	In \Cl{2,0} the general multivector inverse representation is
	\[
	A^{-1}= \frac{1}{a_1^2-a_2^2-a_3^2+a_4^2 }\begin{pmatrix}{a_1} & -{a_2} & -{a_3} & -{a_4}\\
		-{a_2} & {a_1} & -{a_4} & -{a_3}\\
		-{a_3} & {a_4} & {a_1} & {a_2}\\
		{a_4} & -{a_3} & {a_2} & {a_1}\end{pmatrix}\]
	The corresponding characteristic polynomial is
\[	
p(\lambda)={{\left( {{a}_{1}^{2}}-{{a}_{2}^{2}}-{{a}_{3}^{2}}+{{a}_{4}^{2}}-2 {a_1} \lambda+{{\lambda}^{2}}\right) }^{2}}
\]
\end{example}

\begin{example}
In the split quaternions \Cl{1,1} the general inverse representation is
\[
A^{-1}= \frac{1}{-a_1^2+a_2^2-a_3^2+a_4^2}
 \begin{pmatrix}-{a_1} & {a_2} & {a_3} & {a_4}\\
	{a_2} & -{a_1} & {a_4} & {a_3}\\
	-{a_3} & {a_4} & -{a_1} & -{a_2}\\
	{a_4} & -{a_3} & -{a_2} & -{a_1}\end{pmatrix}
\]
	The corresponding characteristic polynomial is
 \[
 	p(\lambda)={{\left( {{a}_{1}^{2}}-{{a}_{2}^{2}}+{{a}_{3}^{2}}-{{a}_{4}^{2}}-2 {a_1} \lambda+{{\lambda}^{2}}\right) }^{2}}
 \]
 
\end{example}

\begin{example}\label{ex:quatinv}
	In the quaternions the general multivector inverse representation is
	\[
	A^{-1} = \frac{1}{a_1^2+a_2^2+a_3^2+a_4^2} \begin{pmatrix}{a_1} & -{a_2} & -{a_3} & -{a_4}\\
		{a_2} & {a_1} & {a_4} & -{a_3}\\
		{a_3} & -{a_4} & {a_1} & {a_2}\\
		{a_4} & {a_3} & -{a_2} & {a_1}\end{pmatrix}
	\]
	corresponding to the characteristic polynomial
\[
	p(\lambda)={{\left( {{a}_{1}^{2}}+{{a}_{2}^{2}}+{{a}_{3}^{2}}+{{a}_{4}^{2}}-2 {a_1} \lambda+{{\lambda}^{2}}\right) }^{2}}
\]

\end{example}
The  FVS algorithm has a direct representation in terms of Clifford multiplications as demonstrated by the next result.
\begin{theorem}[Multivector inversion algorithm]\label{th:clinverse}
	Suppose that $A \in \Cl{p,q}$ is a multivector of maximal grade $r \leq p+q$.
	The Clifford inverse, if it exists, can be computed by the algorithm in $k = 2^{r}$ steps as
	\[
	\begin{array}{l|rr}
		m_1 = A               & t_1 = k \grpart{m_1}_{0},         & c_1= -t_1 \\
		m_2 = A m_1- t_1      & t_2 =\frac{k}{2} {A \ast m_1},    & c_2= -t_2 \\
		$\ldots$ & $ \ldots $ & $ \ldots $ \\
		m_k = A m_{k-1}- t_k  & t_k =  {A \ast m_{k-1}},   & c_k= -t_k
	\end{array}
	\]
	until the step where $m_k=0$ so that
	\begin{flalign}
		 ||A ||= t_k \\ 
		{A}^{-1} = \frac{m_{k-1} }{t_k}  
	\end{flalign}
	The inverse does not exist if $t_k=0$.
	The characteristic polynomial is
	\[
	p(\lambda):= ||\lambda I_k  - A || = \lambda^k + c_1 \lambda^{k-1} + \ldots c_{k-1} \lambda + c_k
	\]
\end{theorem}
\begin{proof}
	The proof follows from the distributivity of the $\pi$ map over the matrix and Clifford products. 
	Furthermore,  $Tr [\mathbf{M}_k]= N \left\langle m_k\right\rangle_0   $.
	Moreover, the FVS algorithm  terminates with $\mathbf{M}_N=0$, which corresponds to $N=2^{p+q}$.	
	
	Suppose that $A$ is of maximal grade $r$. 
	Let $E$ be the set of all blades of grade less or equal to $r$.
	We compute the restricted multiplication tables $\mathbf{M} (E_r)$ and respectively $\mathbf{G} (E_r)$ and form the restricted map $\pi_r$. 
	Then 
	\[
	\pi_r (A A^{-1}) = \pi_r (A)  \pi_r(A^{-1}) =\mathbf{A}  \mathbf{A}^{-1} =\mathbf{I}_n, \quad n=2^r
	\]
	Therefore, the FVS algorithm will terminate in $k=2^r$ steps.  
	Now, suppose that  $t_k \neq 0$, then
	\[
	A m_{k-1}- t_k =0 \Rightarrow A \frac{m_{k-1}}{t_k}=1 \Rightarrow A ^{-1}=	\frac{m_{k-1}}{t_k}
	\]
\end{proof}

\subsection{Low-dimensional inverses}
While inverses can be computed for an arbitrary dimensions the expressions quickly become extremely cumbersome. 
We will consider \Cl{3,0}. Up to sign permutations the results are expected to hold also for \Cl{2,1}, \Cl{1,2}, and \Cl{0,3}.
Let 
\[
A= {a_1}+{e_1} {a_2}+{e_2} {a_3}+{e_3} {a_4}+{a_5} \left( {e_1}\ensuremath{\,}{e_2}\right)  +
{a_6} \left( {e_1}\ensuremath{\,}{e_3}\right) +{a_7} \left( {e_2}\ensuremath{\,}{e_3}\right)
+{a_8} \left( {e_1}\ensuremath{\,}{e_2}\ensuremath{\,}{e_3}\right)
 \] 
Then the application of the LVS algorithm yields 
\[
A^{-1}= \frac{ S + V + BV + PS}{\Delta}
\]
where the determinant is given by
\begin{multline*}
\Delta=
 {{a}_{1}^{4}}-2 {{a}_{1}^{2}}\, {{a}_{2}^{2}}+{{a}_{2}^{4}}-2 {{a}_{1}^{2}}\, {{a}_{3}^{2}}+2 {{a}_{2}^{2}}\, {{a}_{3}^{2}}+{{a}_{3}^{4}} \\
 -2 {{a}_{1}^{2}}\, {{a}_{4}^{2}}+2 {{a}_{2}^{2}}\, {{a}_{4}^{2}}+2 {{a}_{3}^{2}}\, {{a}_{4}^{2}}+{{a}_{4}^{4}}+2 {{a}_{1}^{2}}\, {{a}_{5}^{2}}-2 {{a}_{2}^{2}}\, {{a}_{5}^{2}}-2 {{a}_{3}^{2}}\, {{a}_{5}^{2}}+2 {{a}_{4}^{2}}\, {{a}_{5}^{2}}+{{a}_{5}^{4}}-8 {a_3} {a_4} {a_5} {a_6}+2 {{a}_{1}^{2}}\,  \\
  {{a}_{6}^{2}}-2 {{a}_{2}^{2}}\, {{a}_{6}^{2}}+2 {{a}_{3}^{2}}\, {{a}_{6}^{2}}-2 {{a}_{4}^{2}}\, {{a}_{6}^{2}}+2 {{a}_{5}^{2}}\,  {{a}_{6}^{2}}+{{a}_{6}^{4}}+8 {a_2} {a_4}{a_5} {a_7}-8 {a_2} {a_3} {a_6} {a_7}+2 {{a}_{1}^{2}}\, {{a}_{7}^{2}}+2 {{a}_{2}^{2}}\, {{a}_{7}^{2}}\\
  -2 {{a}_{3}^{2}}\, {{a}_{7}^{2}}-2 {{a}_{4}^{2}}\, {{a}_{7}^{2}}+2 {{a}_{5}^{2}}\, {{a}_{7}^{2}}+2 {{a}_{6}^{2}}\, {{a}_{7}^{2}}+{{a}_{7}^{4}}-8 {a_1} {a_4} {a_5} {a_8}+8 {a_1} {a_3} {a_6} {a_8}-8 {a_1} {a_2} {a_7} {a_8} \\
  +2 {{a}_{1}^{2}}\, {{a}_{8}^{2}}+2 {{a}_{2}^{2}}\, {{a}_{8}^{2}}+2 {{a}_{3}^{2}}\, {{a}_{8}^{2}}+2 {{a}_{4}^{2}}\, {{a}_{8}^{2}}-2 {{a}_{5}^{2}}\, {{a}_{8}^{2}}-2 {{a}_{6}^{2}}\, {{a}_{8}^{2}}-2 {{a}_{7}^{2}}\, {{a}_{8}^{2}}+{{a}_{8}^{4}}
\end{multline*}
and 
\[
S =
 \begin{pmatrix}{{a}_{1}^{3}}-{a_1} {{a}_{2}^{2}}-{a_1} {{a}_{3}^{2}}-{a_1} {{a}_{4}^{2}}+{a_1} {{a}_{5}^{2}}+{a_1} {{a}_{6}^{2}}+{a_1} {{a}_{7}^{2}}-2 {a_4} {a_5} {a_8}+2 {a_3} {a_6} {a_8}-2 {a_2} {a_7} {a_8}+{a_1} {{a}_{8}^{2}}\end{pmatrix}\]
 While the vector part is given by
 \[
 V=
  \begin{pmatrix}
 	e_1 \\
 	e_2 \\
 	e_3 \\
 \end{pmatrix} \begin{pmatrix}-{{a}_{1}^{2}} {a_2}+{{a}_{2}^{3}}+{a_2} {{a}_{3}^{2}}+{a_2} {{a}_{4}^{2}}-{a_2} {{a}_{5}^{2}}-{a_2} {{a}_{6}^{2}}+2 {a_4} {a_5} {a_7}-2 {a_3} {a_6} {a_7}+{a_2} {{a}_{7}^{2}}-2 {a_1} {a_7} {a_8}+{a_2} {{a}_{8}^{2}}\\
 	-{{a}_{1}^{2}} {a_3}+{{a}_{2}^{2}} {a_3}+{{a}_{3}^{3}}+{a_3} {{a}_{4}^{2}}-{a_3} {{a}_{5}^{2}}-2 {a_4} {a_5} {a_6}+{a_3} {{a}_{6}^{2}}-2 {a_2} {a_6} {a_7}-{a_3} {{a}_{7}^{2}}+2 {a_1} {a_6} {a_8}+{a_3} {{a}_{8}^{2}}\\
 	-{{a}_{1}^{2}} {a_4}+{{a}_{2}^{2}} {a_4}+{{a}_{3}^{2}} {a_4}+{{a}_{4}^{3}}+{a_4} {{a}_{5}^{2}}-2 {a_3} {a_5} {a_6}-{a_4} {{a}_{6}^{2}}+2 {a_2} {a_5} {a_7}-{a_4} {{a}_{7}^{2}}-2 {a_1} {a_5} {a_8}+{a_4} {{a}_{8}^{2}}\end{pmatrix}
 \]
 the bi-vector part by 
 \[
 BV =
 \begin{pmatrix}
 	e_1 e_2 \\
 	e_2 e_3 \\
 	e_1 e_3 \\
\end{pmatrix} 
  \begin{pmatrix}-{{a}_{1}^{2}} {a_5}+{{a}_{2}^{2}} {a_5}+{{a}_{3}^{2}} {a_5}-{{a}_{4}^{2}} {a_5}-{{a}_{5}^{3}}+2 {a_3} {a_4} {a_6}-{a_5} {{a}_{6}^{2}}-2 {a_2} {a_4} {a_7}-{a_5} {{a}_{7}^{2}}+2 {a_1} {a_4} {a_8}+{a_5} {{a}_{8}^{2}}\\
 	2 {a_3} {a_4} {a_5}-{{a}_{1}^{2}} {a_6}+{{a}_{2}^{2}} {a_6}-{{a}_{3}^{2}} {a_6}+{{a}_{4}^{2}} {a_6}-{{a}_{5}^{2}} {a_6}-{{a}_{6}^{3}}+2 {a_2} {a_3} {a_7}-{a_6} {{a}_{7}^{2}}-2 {a_1} {a_3} {a_8}+{a_6} {{a}_{8}^{2}}\\
 	-2 {a_2} {a_4} {a_5}+2 {a_2} {a_3} {a_6}-{{a}_{1}^{2}} {a_7}-{{a}_{2}^{2}} {a_7}+{{a}_{3}^{2}} {a_7}+{{a}_{4}^{2}} {a_7}-{{a}_{5}^{2}} {a_7}-{{a}_{6}^{2}} {a_7}-{{a}_{7}^{3}}+2 {a_1} {a_2} {a_8}+{a_7} {{a}_{8}^{2}}\end{pmatrix}
 \]
 and the pseudoscalar part by
\[
PS= e_1 e_2 e_3 \begin{pmatrix}2 {a_1} {a_4} {a_5}-2 {a_1} {a_3} {a_6}+2 {a_1} {a_2} {a_7}-{{a}_{1}^{2}} {a_8}-{{a}_{2}^{2}} {a_8}-{{a}_{3}^{2}} {a_8}-{{a}_{4}^{2}} {a_8}+{{a}_{5}^{2}} {a_8}+{{a}_{6}^{2}} {a_8}+{{a}_{7}^{2}} {a_8}-{{a}_{8}^{3}}\end{pmatrix}
\]
The inverse exists if the determinant $\Delta \neq 0$.
The characteristic polynomial is
\begin{multline*}
p(\lambda)= ({{a}_{1}^{4}}-2 {{a}_{1}^{2}}\, {{a}_{2}^{2}}+{{a}_{2}^{4}}-2 {{a}_{1}^{2}}\, {{a}_{3}^{2}}+2 {{a}_{2}^{2}}\, {{a}_{3}^{2}}+{{a}_{3}^{4}}-2 {{a}_{1}^{2}}\, {{a}_{4}^{2}}+ 
2 {{a}_{2}^{2}}\, {{a}_{4}^{2}}+2 {{a}_{3}^{2}}\, {{a}_{4}^{2}}+{{a}_{4}^{4}}+ \\
 2 {{a}_{1}^{2}}\, {{a}_{5}^{2}}-2 {{a}_{2}^{2}}\, {{a}_{5}^{2}}-2 {{a}_{3}^{2}}\, {{a}_{5}^{2}}+2 {{a}_{4}^{2}}\, {{a}_{5}^{2}}+{{a}_{5}^{4}} 
-8 {a_3} {a_4} {a_5} {a_6}+2{{a}_{1}^{2}}\, {{a}_{6}^{2}}-2 {{a}_{2}^{2}}\, {{a}_{6}^{2}}+2 {{a}_{3}^{2}}\, {{a}_{6}^{2}}-2 {{a}_{4}^{2}}\, {{a}_{6}^{2}}+2 {{a}_{5}^{2}}\, {{a}_{6}^{2}}+ \\
{{a}_{6}^{4}}+ 8 {a_2} {a_4} {a_5} {a_7}-8 {a_2} {a_3} {a_6} {a_7}+2 {{a}_{1}^{2}}\, {{a}_{7}^{2}}+2 {{a}_{2}^{2}}\, {{a}_{7}^{2}}-2 {{a}_{3}^{2}}\, {{a}_{7}^{2}}-2 {{a}_{4}^{2}}\, {{a}_{7}^{2}}+
\\
2 {{a}_{5}^{2}}\, {{a}_{7}^{2}}+2 {{a}_{6}^{2}}\, {{a}_{7}^{2}}+{{a}_{7}^{4}}-8{a_1} {a_4} {a_5} {a_8}+8 {a_1} {a_3} {a_6} {a_8}-8 {a_1} {a_2} {a_7} {a_8}+2 {{a}_{1}^{2}}\, {{a}_{8}^{2}}+2 {{a}_{2}^{2}}\, {{a}_{8}^{2}}+2 {{a}_{3}^{2}}\, {{a}_{8}^{2}}+
\\
2 {{a}_{4}^{2}}\, {{a}_{8}^{2}}-2 {{a}_{5}^{2}}\, {{a}_{8}^{2}}-2 {{a}_{6}^{2}}\, {{a}_{8}^{2}}-2 {{a}_{7}^{2}}\, {{a}_{8}^{2}}+{{a}_{8}^{4}}+ \\
( -4 {{a}_{1}^{3}}+4 {a_1} {{a}_{2}^{2}}+4 {a_1} {{a}_{3}^{2}}+4 {a_1} {{a}_{4}^{2}}-4 {a_1} {{a}_{5}^{2}}-4 {a_1} {{a}_{6}^{2}}-4 {a_1} {{a}_{7}^{2}}+ 
8 {a_4} {a_5} {a_8}-8 {a_3} {a_6} {a_8}+8 {a_2} {a_7} {a_8}-4 {a_1} {{a}_{8}^{2}}) \lambda+ \\
\left( 6 {{a}_{1}^{2}}-2 {{a}_{2}^{2}}-2 {{a}_{3}^{2}}-2 {{a}_{4}^{2}}+2 {{a}_{5}^{2}}+2 {{a}_{6}^{2}}+2 {{a}_{7}^{2}}+2 {{a}_{8}^{2}}\right)  {{{{\lambda}}}^{2}}-4 {a_1} {{{{\lambda}}}^{3}}+{{{{\lambda}}}^{4}})^2
\end{multline*} 

\subsection{Numerical examples} 
\begin{example}
	Let us compute an example in \Cl{3,1}, having rational coefficients.
	Let 
	\begin{multline*}
	A= -6+\frac{1}{5}\, {e_1}-{e_2}+3 {e_3}-\frac{3}{5}\, {e_4}-{e_1}\,{e_2}-7 \left( {e_1}\,{e_2}\,{e_3}\right) -\frac{1}{3}\, \left( {e_1}\,{e_2}\,{e_3}\,{e_4}\right) + \\
	{e_1}\,{e_2}\,{e_4}+\frac{3}{2}\, \left( {e_1}\,{e_3}\right) +2 \left( {e_1}\,{e_3}\,{e_4}\right) +3 \left( {e_1}\,{e_4}\right) +{e_2}\,{e_3}+\frac{7}{6}\, \left( {e_2}\,{e_3}\,{e_4}\right) -{e_2}\,{e_4}+\frac{7}{2}\, \left( {e_3}\,{e_4}\right)
	\end{multline*}
	The computed inverse is 
	\begin{multline*}
	A^{-1}=-\frac{1188400}{12512583}-\frac{618140}{37537749}\, {e_1}+\frac{3103100}{112613247}\, {e_2}+\frac{709300}{12512583}\, {e_3} \\
	-\frac{1187740}{112613247}\, {e_4}+\frac{4015700}{112613247}\, \left( {e_1}\,{e_2}\right) +\frac{1458400}{12512583}\, \left( {e_1}\,{e_2}\,{e_3}\right) +\frac{1741900}{37537749}\, \left( {e_1}\,{e_2}\,{e_3}\,{e_4}\right) \\
	-\frac{159700}{4170861}\, \left( {e_1}\,{e_2}\,{e_4}\right) +\frac{434450}{12512583}\, \left( {e_1}\,{e_3}\right) -\frac{597200}{12512583}\, \left( {e_1}\,{e_3}\,{e_4}\right) -\frac{6667300}{112613247}\, \left( {e_1}\,{e_4}\right)\\
	 +\frac{1366700}{37537749}\, \left( {e_2}\,{e_3}\right) 
	 +\frac{224950}{37537749}\, \left( {e_2}\,{e_3}\,{e_4}\right) +\frac{280100}{37537749}\, \left( {e_2}\,{e_4}\right) -\frac{1817150}{37537749}\, \left( {e_3}\,{e_4}\right) 
	\end{multline*}
   The intermediate steps of the computations are available for inspection by setting the appropriate flags in the code.
   However, the rational numbers involved quickly become very long.
   
\end{example}

Let us compute several examples in \Cl{5,2}.
\begin{example}
 Let 
 \[
 A= 1 -{e_2} +  {e_1} \, {e_2}\,{e_3}\,{e_4} \, {e_5}\,{e_6}\,{e_7}
 \]
 Then
 \[
 A^{-1}= \frac{1}{5}-\frac{1}{5}\, {e_2}-\frac{3}{5}\, \left( {e_1}\,{e_2}\,{e_3}\,{e_4}\,{e_5}\,{e_6}\,{e_7}\right) +\frac{2}{5}\, \left( {e_1}\,{e_3}\,{e_4}\,{e_5}\,{e_6}\,{e_7}\right) 
 \] 
\end{example}

\begin{example}
	Let 
	\[
	A= 1 -{e_2} -\frac{1}{8} e_1\ e_5 +  {e_1} \, {e_2}\,{e_3}\,{e_4} \, {e_5}\,{e_6}\,{e_7}
	\]
	Then
	\begin{multline*}
		A^{-1}=
		\frac{20544}{102785}-\frac{20544}{102785} {e_2}-\frac{3932224}{6475455}   {e_1} {e_2} {e_3} {e_4} {e_5} {e_6} {e_7}  +\frac{459776}{6475455}   {e_1} {e_2} {e_5}  +\frac{2646016}{6475455}  {e_1} {e_3} {e_4} {e_5} {e_6} {e_7}  \\
		-\frac{362504}{6475455}   {e_1} {e_5} -\frac{1024}{102785}  {e_2} {e_3} {e_4} {e_6} {e_7}  +\frac{1024}{102785}  {e_3} {e_4} {e_6} {e_7} 
	\end{multline*}
	Let 
	\[
	A= 1-e_2+ e_1 \, e_2 \, e_3\, e_4 \, e_5 \, e_6 \ e_7 - e_1 \, e_4 \, e_5
	\]
	Then
	\[
	A^{-1}=
	\frac{1}{2}(1-e_1 \, e_2 \, e_3 \, e_4 \, e_5 \, e_6 \, e_7+e_1 \, e_3 \, e_4 \, e_5 \, e_6 \, e_7+e_2 \, e_3 \, e_6 \, e_7)
	\]
\end{example}

\begin{example}
In \Cl{2,2}let 
\[
A= 1+{e_1}+{e_1}\,{e_3}\,{e_4}-2 \left( {e_2}\,{e_3}\right) 
\]
\begin{flalign*}
	t_1=-8, &\quad m_1=1+{e_1}+{e_1}\,{e_3}\,{e_4}-2 \left( {e_2}\,{e_3}\right)   \\
	t_2=12, &\quad m_2=-3-6 {e_1}-4 \left( {e_1}\,{e_2}\,{e_3}\right) -6 \left( {e_1}\,{e_3}\,{e_4}\right) +12 \left( {e_2}\,{e_3}\right) +2 \left( {e_3}\,{e_4}\right) \\
	t_3=40, &\quad m_3=-15+9 {e_1}+20 \left( {e_1}\,{e_2}\,{e_3}\right) +5 \left( {e_1}\,{e_3}\,{e_4}\right) -10 \left( {e_2}\,{e_3}\right) -10 \left( {e_3}\,{e_4}\right) \\
	t_4=-98, &\quad m_4=49+4 {e_1}-8 \left( {e_1}\,{e_2}\,{e_3}\right) +20 \left( {e_1}\,{e_3}\,{e_4}\right) -40 \left( {e_2}\,{e_3}\right) +4 \left( {e_3}\,{e_4}\right)\\
	t_5=-24, &\quad m_5=15-33 {e_1}-56 \left( {e_1}\,{e_2}\,{e_3}\right) -25 \left( {e_1}\,{e_3}\,{e_4}\right) +50 \left( {e_2}\,{e_3}\right) +28 \left( {e_3}\,{e_4}\right)\\
	t_6=156, &\quad m_6=-117+42 {e_1}+60 \left( {e_1}\,{e_2}\,{e_3}\right) -6 \left( {e_1}\,{e_3}\,{e_4}\right) +12 \left( {e_2}\,{e_3}\right) -30 \left( {e_3}\,{e_4}\right) \\
	t_7=-72, &\quad m_7=63-9 {e_1}-12 \left( {e_1}\,{e_2}\,{e_3}\right) +3 \left( {e_1}\,{e_3}\,{e_4}\right) -6 \left( {e_2}\,{e_3}\right) +6 \left( {e_3}\,{e_4}\right) 
\end{flalign*}
so that 
\[
A^{-1}=-\frac{1}{3} \left( -3-3 {e_1}-4 \left( {e_1}\,{e_2}\,{e_3}\right) +{e_1}\,{e_3}\,{e_4}-2 \left( {e_2}\,{e_3}\right) +2 \left( {e_3}\,{e_4}\right) \right) \]

On the other hand, for certain low-dimensional cases -- n=2 and n=4, the FLV algorithm can be simplified as follows.
\begin{flalign*}
	t_1=4, &\quad m_1=1+{e_1}+{e_1}\,{e_3}\,{e_4}-2 \left( {e_2}\,{e_3}\right)  \\
	t_2=2, &\quad m_2=1-2 {e_1}-4 \left( {e_1}\,{e_2}\,{e_3}\right) -2 \left( {e_1}\,{e_3}\,{e_4}\right) +4 \left( {e_2}\,{e_3}\right) +2 \left( {e_3}\,{e_4}\right) \\
	t_3=-12, &\quad m_3=-9+3 {e_1}+4 \left( {e_1}\,{e_2}\,{e_3}\right) -{e_1}\,{e_3}\,{e_4}+2 \left( {e_2}\,{e_3}\right) -2 \left( {e_3}\,{e_4}\right) 
\end{flalign*}
leading to the same result for $A^{-1}$.	
\end{example}

\section{Discussion}\label{sec:disc}

The main theoretical objective of the present paper was to demonstrate a universal construction of a matrix algebra representation of 
Clifford algebra of non-degenerate arbitrary signature \textit{(p, q)}.
This is achieved by exhibiting an explicit isomorphism between a given Clifford algebra and its faithful matrix representation.
Instrumental for the presented approach is the Clifford algebra matrix multiplication table  $\mathbf{M}$.
The structure of this table is  constrained by the structure of the Clifford algebra itself, 
which allows for deriving useful relations between the matrix entries thus reducing the computational complexity of the matrix construction. 
As  demonstrated further in the paper, the matrix multiplication table $\mathbf{M}$ encodes all  properties of the algebra
It should be noted that such a representation is not optimal from data compression perspective.  
On the other hand, it can be used for automatic computer code generation for the lower-dimensional algebras. 
Target applications of such an approach can be computational environments, such as Matlab, or computer code preprocessors, such as  Gaalop for C++ \cite{Hildenbrand2010}.  

It should be noted that the FLV algorithm  can result easily in float precision issues.  
On the other hand, the FLV algorithm is in fact a proof certificate, 
that is it terminates in a finite number of steps and if an inverse exits it provides it.

\section*{Funding} 
Funding for the initial \textit{clifford} development was provided by a grant from Research Fund -- Flanders (FWO), contract numbers G.0C75.13N, VS.097.16N.  

\section*{Conflicts of interest}
The author declares no competing intesrests. 

\section*{Availability of data and material}
Data are available in a Zenodo repository \cite{Prodanov2016}.

\section*{Code availability} 
Maxima code is available in a Zenodo repository \cite{Prodanov2016}.


\appendix

\section{Auxiliary results}\label{sec:properties}

\begin{proposition}[Maximal element]\label{prop:pscalar}
	The algebra \Cl{p,q,r}  has a maximal element $	\mathcal{I} = e_1 \ldots e_n$ called pseudoscalar. 
\end{proposition}
\begin{proof}
	Suppose that $r>0$. Then if a product contains more than one nilpotent generators of the same index the simplified form is 0 by Lemma \ref{th:permeq}. 
	If a product contains exactly one nilpotent generator per index the maximal element is the product of all generators by  Th. \ref{th:simpf}.
	\[
	\mathcal{I} = e_1 \ldots e_n, \ n=p+q+r
	\]
	Suppose that $r=0$. Then by Th. \ref{th:simpf} the maximal element is the product of all generators
	\[
	\mathcal{I} = e_1 \ldots e_n, \ n=p+q
	\]
	This element is referred to as the \textit{pseudoscalar} of the algebra.
\end{proof}


\begin{proposition}\label{prop:antisym}
	For generator elements $e_s$ and $e_t$
	\[
	\mathbf{E}_s \mathbf{E}_t + \mathbf{E}_t \mathbf{E}_s= \mathbf{0}
	\]
\end{proposition}
\begin{proof}
	Consider the basis elements $e_s$ and $e_t$. 
	By linearity and Clifford-product distributivity of the $\pi$ map
	\[
	\pi: e_s e_t + e_t e_s = 0 \mapsto \pi(e_s e_t) + \pi(e_t e_s) = \mathbf{0}
	\]
	Therefore, for vector elements
	\[
	\mathbf{E}_s \mathbf{E}_t + \mathbf{E}_t \mathbf{E}_s= \mathbf{0}
	\]
\end{proof}

\begin{proposition}\label{prop:vsquare}
	$ \mathbf{E_s} \mathbf{E_s} = \sigma_s \mathbf{I} $ 
\end{proposition}
\begin{proof}
	Consider the matrix
	$ \mathbf{W} = \mathbf{G}  \mathbf{A_s} \mathbf{G}  \mathbf{A_s}   $.
	Then element-wise
	$w_{\mu \nu} = \sum_{\lambda} \sigma_\mu \sigma_\lambda a_{ \mu \lambda } a_{  \lambda \nu}$.
	By Lemma \ref{th:sparsity} $\mathbf{W} $ is sparse so that
	$w_{\mu \nu} =  (0; \sigma_\mu \sigma_q a_{ \mu q } a_{  q \nu}) $.
	
	From the structure of $\mathbf{M}$ for the entries containing the element $e_S$ we have the equivalence
	\[
	\begin{cases}
		e_{M} e_{Q} = a^s_{\mu q} e_S,  \ \  S= M \triangle Q \\
		e_{Q} e_{M} = a^s_{q \mu} e_S, 
	\end{cases}	
	\]
	After multiplication of the equations we get
	\[
	e_{M} e_{Q}  e_{Q} e_{M} = a^s_{\mu q} e_S a^s_{q \mu} e_S 
	\] 
	which simplifies to the \textit{First fundamental identity}:
	\begin{equation}\label{eq:fistId}
		\sigma_q \sigma_\mu = a^s_{\mu q}   a^s_{q \mu} \sigma_s
	\end{equation}
	We observe that if  $\sigma_\mu=0$ or $\sigma_q =0$ the result follows trivially.
	In this case also $\sigma_s=0$.
	
	Therefore, let's suppose that $ \sigma_s \sigma_q \sigma_\mu \neq 0 $.
	We multiply both sides by $ \sigma_s \sigma_q \sigma_\mu  $
	\[
	\sigma_s = \sigma_q \sigma_\mu a^s_{\mu q}   a^s_{q \mu} 
	\]
	However, the RHS is a diagonal element of $\mathbf{W}$, therefore by the sparsity it is the only non-zero element for a given row/column so that 
	\[
	\mathbf{W}= \mathbf{E_s^2} = \sigma_s \mathbf{I}
	\] 
\end{proof}

\section{Properties of the scalar product multiplication tables}\label{sec:scalartab}

\begin{theorem}
	The scalar product table is a diagonal matrix $\mathbf{G}$, which is invariant under orthogonal transformations.
\end{theorem}
\begin{proof}
	The proof is based on Macdonald \cite{Macdonald2002}.
	From the definition of the scalar product it is obvious that $\mathbf{G}$ is diagonal. 
	Consider the orthogonal transformation $ \mathbf{F} = \mathbf{A} \mathbf{E}$ with orthogonal matrix $\mathbf{A}$.
	We evaluate $\mathbf{F}^T \mathbf{F}$.
	Then element-wise (summation by repeated indices)
	\begin{equation}
	f_n f_m = a_{n j} a_{j m} e_N  e_M
	\end{equation}
	Then for $m=n$
	\[
	f_n f_n = a_{n j} a_{j n} e_N  e_N = \sigma_n
	\]
	by the orthogonality of entries. 
	Then $
	diag \left\langle  \mathbf{F}^T \mathbf{F} \right\rangle = \mathbf{G}
	$, where \textit{diag} denotes the diagonal elements.
\end{proof}

\begin{definition}[Sparsity property]
	A matrix has the \textit{sparsity property} if it has exactly one non-zero element per column and exactly one non-zero element per row. Such a matrix we call sparse.
\end{definition}

\begin{proposition}\label{th:sparsec}
	For $e_s \in \mathbf{E}$ the matrix $\mathbf{A}_s = C_{e_s} (\mathbf{M})$ is sparse. 
\end{proposition}
\begin{proof}
	Fix an element $e_s \in \mathbf{E}$.
	Consider a row $k$.
	By prop \ref{prop:msparse} there is a $j$, such $e_{k j}= e_s$.
	Then $a_{k j}= \mu_{k j}$, while for $i \neq j$ $a_{k i}=0$.
	Consider a column $m$
	By prop \ref{prop:msparse} there is a $j$, such $e_{j m}= e_s$.
	Then $a_{j m}= \mu_{j m}$, while for $i \neq j$ $a_{i m}= 0$.
	Therefore, $\mathbf{A}_s$ has the sparsity property.
\end{proof}

\begin{lemma}[Sparsity lemma]\label{th:sparsity}
	If the matrices $\mathbf{A}$ and $\mathbf{B}$ are sparse then so is $\mathbf{C}=\mathbf{A B}$.
	Moreover,
	\[
	c_{ij} = \begin{cases}
	0 \\
	a_{i q} b_{q j} 
	\end{cases}
	\] 
	(no summation!) for some index $q$.
\end{lemma}
\begin{proof}
	Consider two sparse square matrices $\mathbf{A}$ and $\mathbf{B}$  of dimension $n$.
	Let $c_{ij} = \sum_{\mu} a_{i \mu} b_{\mu j}$. Then
	 as we vary the row  index $i$ then there is only one index $ q \leq n$, such that $ a_{i q}  \neq 0$.
	As we vary the column index $j$ then there is only one index  $ q \leq n$, such that $ b_{ q j}  \neq 0$.
	Therefore, $c_{ij} = (0; a_{i q} b_{q j}  ) $
	for some $q$ by the sparsity of $\mathbf{A}$ and $\mathbf{B}$.

	As we vary the row  index $i$ then $c_{q j} =0 $ for $i \neq q$ for the column $j$ by the sparsity of  $\mathbf{A}$.
	As we vary the column index $j$ then $c_{i q}=0$ for $j \neq q$ for the row $i$ by the sparsity of  $\mathbf{B}$.
	Therefore, $\mathbf{A B}$ is sparse. 
\end{proof}
\begin{corollary}\label{cor:sp1}
	If $\mathbf{A}$ is sparse then $\mathbf{A^2}$ is diagonal.
\end{corollary}
\begin{proof}
	Let $\mathbf{C}=\mathbf{A^2}$
	We observe that by sparsity $a_{i q} b_{q j} \neq 0$ if $i=j$.
\end{proof}

\subsection{Non-degenerate algebras}
\label{sec:nondegen}

\begin{proposition}
	For non degenerate algebras \Cl{p,q} the matrix $\mathbf{G}$ is orthogonal and 	
	\[
	\mathbf{G} \mathbf{G} = \mathbf{I}
	\]
\end{proposition}
\begin{proof}
	$\mathbf{G} = \mathbf{G} ^T $ by symmetry. Since its elements are valued $\pm 1$  then element-wise
	$
	\sum_{\lambda} g_{\mu \lambda} g_{\lambda \mu} = \sigma_\mu^2 = 1
	$
\end{proof}

\begin{proposition}\label{prop:einv} 
	For non-degenerate algebras \Cl{p, q}
	$
	\mathbf{E_s^{-1}} = \sigma_s \mathbf{E_s}
	$.
\end{proposition}

\begin{theorem}
	For non-degenerate algebras \Cl{p,q} 
	the matrices $\mathbf{A}_s$ and  $\mathbf{E}_s$ are orthogonal.
\end{theorem}
\begin{proof}
	First we prove an auxiliary equation.
	Consider the matrix
	\[
     \mathbf{Y} =\mathbf{ G A}_s - \sigma_s \mathbf{ A}_s^T \mathbf{G }
	\]
	Element-wise
	$ y_{\mu \nu} = \sigma_\mu a^s_{\mu \nu}  - \sigma_s   a^s_{\nu \mu} \sigma_\nu $. 
	Suppose that $a^s_{\nu \mu} \neq 0 $.
	Then $$ 
	y_{\mu \nu} = a^s_{\nu \mu} \left( \sigma_\mu a^s_{\mu \nu} a^s_{\nu \mu}  - \sigma_s \sigma_\nu    \right)
	= \sigma_\nu  a^s_{\nu \mu}  \underbrace{\left(\sigma_\mu a^s_{\mu \nu} a^s_{\nu \mu} \sigma_\nu   - \sigma_s \right)}_{ 0}  = 0 
	 $$
	 Therefore,
	 \begin{equation}\label{eq:Greverese}
	  \mathbf{ G A}_s = \sigma_s \mathbf{ A}_s^T \mathbf{G }
	 \end{equation}
	 Right multiply by $ \mathbf{ G A}_s $:
	 \[
	 \underbrace{\mathbf{ G A}_s \mathbf{ G A}_s}_{\mathbf{E^2}_s} = \sigma_s \mathbf{ A}_s^T \underbrace{\mathbf{G } \mathbf{ G} }_{\mathbf{I}} \mathbf{ A}_s
	 \]
	 so that
	 \[
	 \sigma_s \mathbf{I} = \sigma_s \mathbf{ A}_s^T \mathbf{ A}_s
\Leftrightarrow
	 \mathbf{ A}_s^T \mathbf{ A}_s = \mathbf{I}
	 \]
	 Since $\mathbf{G }$ is orthogonal the second assertion follows as well. 
\end{proof}
\begin{corollary}[Second fundamental identity]
	Consider the element $e_s$ in a non-degenerate algebra \Cl{p,q}. Then there is an index 
	$q \leq n= p+q $, such that
	\[
	\begin{cases}
		a^s_{\mu q} a^s_{q \mu} & =  \sigma_s \\
		\sigma_\mu \sigma_q & = 1	
	\end{cases}
	\]
\end{corollary}
\begin{proof}
 $ \mathbf{A}_s {\mathbf{A}_s^T} = \mathbf{I} \Rightarrow (a^s_{q \mu})^2 =1 $ for some $q$.
 $ \mathbf{E}_s {\mathbf{E}_s^T} = \mathbf{I} \Rightarrow \sigma_\mu \sigma_q (a^s_{q \mu})^2 = \sigma_\mu \sigma_q = 1 $ for some $q$. The final assertion follows from the first fundamental identity.
\end{proof}

\subsection{Degenerate algebras}
\label{sec:degen}
For degenerate algebras a weaker version of the result holds
\begin{proposition}
	Consider a degenerate algebra \Cl{p,q,r} and let $\mathbf{H}:=\mathbf{G}^2$. Then
	$\mathbf{H}$ is idempotent. 	
\end{proposition}
\begin{proof}
	We observe that since $\mathbf{G}$ is diagonal and $\sigma_i ^2=1$ it follows that
	$\mathbf{H}$ contains only 1's and 0's on its diagonal
	so $h_{\mu \lambda} h_{\lambda \mu} = h_{\mu \mu}^2 \in \{0, 1 \}$  and
	\[
	\mathbf{H} \mathbf{H} = \mathbf{H}
	\]
\end{proof}


\section{Maxima code}
\label{sec:code}

The code was implemented using the Maxima package \symb{clifford} developed by the author \cite{Prodanov2016,Prodanov2016a}. 
Examples presented above can be recomputed using the command \symb{elem2mat1} after proper initialization of the package. 
The key functions of the implementation are listed below.

\begin{lstlisting}[caption={Algebra construction code}, label=lst:rep]

/*  computes blade representation */
climatrep(vv):=block([n, q, AA, lst: elements(all), G, EE],
	local(AA, G, EE),
	n:length(lst),
	/* multiplication table of the algebra */
	AA:genmatrix( lambda([i,j], dotsimpc(lst[i] . lst[j] ) ), n),
	/* signature of the algebra */
	G:diag(AA),
	EE:matrixmap(lambda([q],  scalarpart(ratcoeff (vv, q))), AA),
	G.EE
);

/*  computes expression representation */
elem2mat1(expr, [ulst] ):=block(
	if emptyp(ulst) then ulst:false 
	else ulst:true,
	if mapatom(expr) then return(climatrep1(expr)),
	if ulst=true then
		maplist(climatrep1, expr)
	else
		map(climatrep1, expr)
);

 
/*  Faddeev Le Verrier algorithm */
fadlevi(A):=block([n, m, c, M, K, I, i, cq],
	if not matrixp(A) then error ("not a matrix ", A),
	[n,m]:matrix_size(A),
	if n#m then error ("not a square matrix"),
	c: makelist(1,n+1),
	M:ident(n),
	I:ident(n),
	for i:1 thru n do (
		K: ratsimp(A.M),
		c [i+1]: -1/i * mattrace(K),
		if i<n then  
			M:ratsimp(K + c[i+1]*I)
	), 
	c:factor(expand(c)),
	cq:-last(c),
	if cq=0 then cq:1,
	M: ratsimp(M/cq),
	[M, c]
);

/* pseudoinverse : for n idempotent it returns the complement zero divisor */
fadlevicg1(A):=block([ M:1, K, i:1, n, k: maxgrade(A), cq],
	n: 2^k,
	if _debug1=true then print("n=",n),
	A:rat(A),
	while i<n and K#0 do (
		K: dotsimpc(expand (A.M)),
		cq: -n/i* scalarpart(K),
		if _debug1=true then print(i, " ", cq),
		if _debug1=all then print("t[",i,"]=", cq, " m[",i,"]=", K),
		if K#0 then
			M: rat(K + cq), 
		i:i+1
	), 
	K: dotsimpc(expand(A.M)),
	cq: -n/i* scalarpart(K),
	cq:factor(cq),
	if cq=0 then cq:1,
	M:factor(-(M)/cq)
);
\end{lstlisting}

\bibliographystyle{spmpsci}
\bibliography{clibib1}

\end{document}